\documentclass{tran-l}

\usepackage{amssymb}

\usepackage{graphicx}

\usepackage[cmtip,all]{xy}

\usepackage{amsmath}
\usepackage{mathrsfs}
\usepackage[colorlinks=true, urlcolor=blue,bookmarks=true,bookmarksopen=true, citecolor=blue,hypertex]{hyperref}

\newtheorem{theorem}{Theorem}[section]
\newtheorem{lemma}[theorem]{Lemma}
\newtheorem{corollary}[theorem]{Corollary}
\newtheorem{prop}[theorem]{Proposition}
\newtheorem{propdef}[theorem]{Proposition (Definition)}

\theoremstyle{definition}
\newtheorem{definition}[theorem]{Definition}

\newtheorem{assumption}[theorem]{Assumption}
\newtheorem{notation}{Notation}[section]

\theoremstyle{remark}
\newtheorem{remark}[theorem]{Remark}

\numberwithin{equation}{section}

\newcommand{\idelt}{1\!\!1}
\begin{document}

\title{A relative Seidel morphism and the Albers map}


\author{Shengda Hu}
\address{Department of pure mathematics, University of Waterloo, Waterloo (Ontario), Canada}
\curraddr{}
\email{hshengda@math.uwaterloo.ca}
\thanks{}

\author{Fran\c cois Lalonde}
\address{D\'epartement de math\'ematiques et de statistique, Universit\'e de Montr\'eal, Montr\'eal (Qu\'ebec), Canada}
\curraddr{}
\email{lalonde@dms.umontreal.ca}
\thanks{}

\subjclass[2000]{53D12, 53D40, 53D45, 57R58, 57S05}

\date{}

\dedicatory{}

\begin{abstract} In this note, we introduce a relative (or Lagrangian) version of the Seidel homomorphism that assigns to each homotopy class of paths in ${\rm Ham}(M)$, starting at the identity and ending on the subgroup that preserves a given Lagrangian submanifold $L$, an element in the Floer homology of $L$. We show that these elements are related to the absolute Seidel elements by the Albers map. We also study, for later use, the effect of reversing the signs of the symplectic structure as well as the orientations of the generators and of the operations on the Floer homologies.
\end{abstract}

\maketitle
\tableofcontents
%

\section{Introduction}\label{intro}

Let $(M, \omega)$ be a closed symplectic manifold. Seidel constructed in \cite{Seidel} a map from a covering of $\pi_1{\rm Ham}(M, \omega)$ to the invertible elements of either $QH_*(M, \omega)$ or $FH_*(M)$. It has found many uses, for example, in the study of Hamiltonian fibrations (Lalonde--McDuff--Polterovich  \cite{LalondeMcDuffPolterovich}) and the quantum ring structure of toric varieties (McDuff-Tolman \cite{McDuffTolman}). In this article, we introduce a similar construction when given a Lagrangian submanifold $L$ in $(M, \omega)$. Instead of considering the loops in ${\rm Ham}(M, \omega)$, we consider the paths in ${\rm Ham}(M, \omega)$ starting at the identity and ending in the subgroup ${\rm Ham}_L(M, \omega)$:
$${\rm Ham}_L(M, \omega) := \{\varphi \in {\rm Ham}(M, \omega) | \varphi (L) = L \}$$
Another natural subgroup to consider consists of Hamiltonian symplectomorphisms that fix $L$ pointwisely. It is easy to see that any diffeomorphism of $L$ that is isotopic to the identity can be extended to a Hamiltonian symplectomorphism in $M$ that preserves $L$. It follows that the two choices of subgroups essentially differ by ${\rm Diff}(L)$. For the purpose of this paper, the reader can think of either choice. 

Under the monotonicity assumption in the Lagrangian setting, one can define the Seidel elements for the elements in a covering of $\pi_1({\rm Ham}(M, \omega), {\rm Ham}_L(M, \omega))$. There is a homotopy exact sequence for the Hamiltonian groups, and we show that
the following diagram commute (cf. corollary \ref{lagSeidel:3actions})
\begin{equation}\label{intro:commutdiag}
\text{
\xymatrix{
\widetilde \pi_1{\rm Ham}(M, \omega) \ar[r]\ar[d]_{\Psi} & \widetilde \pi_1({\rm Ham}(M, \omega), {\rm Ham}_L(M, \omega)) \ar[d]_{\Psi_L} 
\ar[r] & \widetilde \pi_0{\rm Ham}_L(M, \omega)\\
FH_*(M) \ar[r]^{\mathscr A} & FH_*(M, L) 
&
}}
\end{equation}
where $\Psi$ and $\Psi_L$ denotes the respective absolute and relative Seidel maps and $\mathscr A$ denotes Albers' comparison map between $FH_*(M)$ and $FH_*(M, L)$ \cite{Albers}.

We should explain the above diagram a little more. The Seidel maps are defined for the extensions $\widetilde\pi_1$ of the respective $\pi_1$'s by the corresponding period groups $\Gamma_\omega$ or $\Gamma_L$. 
An element $\widetilde g \in \widetilde \pi_1{\rm Ham}(M, \omega)$ can be viewed canonically as an element in $\widetilde \pi_1({\rm Ham}(M, \omega), {\rm Ham}_L(M, \omega))$ (cf. lemma \ref{lagSeidel:groupinclusion}) and the corresponding Seidel elements are related by $\mathscr A$. 
The group $\widetilde \pi_0{\rm Ham}_L(M, \omega)$ is an extension of $\pi_0 {\rm Ham}_L(M, \omega)$ 
so that the top sequence is exact. 
On the other hand, we may adopt McDuff's point of view in \cite{McDuff} where the Seidel map is defined on $\pi_1$'s directly by choosing a prefered extension $\widetilde g$ for each $g \in \pi_1$. Then \eqref{intro:commutdiag} holds with the non-extended homotopy groups.

One of the original motivations for studying these Seidel elements and the above diagram was to obtain information on the third term in the exact sequence, namely $\pi_0{\rm Ham}_L(M, \omega)$. This is the most elementary question that one may ask about the subgroup ${\rm Ham}_L(M, \omega)$ of Hamiltonian transformations leaving a given Lagrangian submanifold invariant. 

If, say, one calls the elements in the image of ${\mathscr A} \circ \Psi$ the {\it circular Seidel elements} in $FH_*(M,L)$ and the elements in the image of $\Psi_L$ the {\it semi-circular} ones, then by the commutativity of the diagram, the latter ones contain the former ones. One could try to find semi-circular, but not circular, elements by computing explicitly the two Seidel morphisms. This would imply that 
$$
\pi_1{\rm Ham}(M, \omega) \rightarrow \pi_1({\rm Ham}(M, \omega), {\rm Ham}_L(M, \omega))
$$
is not onto, and therefore $\pi_0{\rm Ham}_L(M, \omega)$ is not trivial. Of course, if a given component of ${\rm Ham}_L(M, \omega)$ is made of Hamiltonian diffeomorphisms whose restrictions to $L$ is not isotopic to the identity, it is obviously not the identity component. Hence our construction is useful  when  components are made of Hamiltonian diffeomorphisms whose restriction to $L$ is isotopic to the identity (or cannot be easily shown to be non-isotopic to the identity).

In presence of a real structure on $M$, i.e. of an anti-symplectic
involution $c$ with fixed point set $L$,
one could replace everywhere ${\rm Ham}_L(M, \omega)$ by its subgroup
${\rm Ham}_c(M, \omega)$ of
Hamiltonian diffeomorphisms of $M$ commuting with $c$. Obviously the
corresponding diagram \eqref{intro:commutdiag} would still commute. In this paper, we restrict ourselves to a detailed setting of the theory, postponing to a forthcoming paper applications and computations of examples.

In section \S\ref{lagFloer}, we review Lagrangian Floer homology in the setting of Hamiltonian paths, cf. \cite{Albers}. Other versions of Floer homology of a single Lagrangian already exist, e.g. \cite{Oh, BiranCornea, BiranCornea2, FOOO} by Oh, Fukaya--Oh--Ohta--Ono, Biran--Cornea. Because we will be working over ${\mathbb R}$, we include a discussion of the coherent orientation for the Floer trajectories, essentially following Fukaya--Oh--Ohta--Ono \cite{FOOO}. The necessary gluing is similar to those found in Albers \cite{Albers}. The half-pair-of-pants product is analogous to the pair-of-pants product in the Hamiltonian Floer homology with some notorious differences (say its non-commutativity). It should coincide with the product defined from holomorphic triangles as in \cite{FOOO} or using the linear cluster complex as in \cite{BiranCornea}. For later use, we also discuss the action of $FH_*(M)$ as well as the Albers' comparison map. Other versions of the action of $FH_*(M)$ has been described before, e.g. in the context of the linear cluster complex (or pearl complex) \cite{BiranCornea} or in the context of $L_\infty$-action on the $A_\infty$-algebra in \cite{FOOO}.

In \S\ref{lagSeidel} we carry out Seidel's construction for $FH_*(M, L)$ and show that it has the expected properties. It is also in this section that we show the commutativity of diagram \eqref{intro:commutdiag}.   Finally, for later use related to the computations in $(X, \omega) \times (X, -\omega)$, we explain in \S\ref{reversed:novikov} and in  \S\ref{reversed:lagFloer} the effect of reversing the orientations of the generators of the symplectic and Lagrangian Floer  homologies a well as the reversal of time in operations on Floer homologies. 

   We would like to mention here that the possibility of defining a relative Seidel morphism appears implicitly  in the recent paper of R\'emi Leclercq in \cite{Leclercq}. Indeed, the proof of his basic proposition 3.1 that he needs to define his Lagrangian spectral invariants, contains the main ingredients of the construction of a relative Seidel morphism, even though it is not presented in these terms. 

\vspace{0.1in}
\noindent
{\bf Acknowledgement.} We would like to thank the referee for the constructive suggestions. The first author would like to thank Octav Cornea, Dusa McDuff and Jean-Yves Welschinger for valuable discussions.

\section{Lagrangian Floer theory}\label{lagFloer}

Let $(M, \omega)$ be a symplectic manifold and $L$ a Lagrangian submanifold. We set up here the Floer theory for $(M, L)$ with a generic Hamiltonian purturbation. 

\subsection{Novikov rings}\label{lagFloer:novikov}
We think of $\omega$, $c_1(TM)$ and $\mu_L$ as functions on $\pi_2(M)$ or $\pi_2(M, L)$ and denote them as:
$$I_\omega : \pi_2(M) \text{ or } \pi_2(M, L) \to {\mathbb R}, I_c : \pi_2(M) \to {\mathbb R} \text{ and } I_\mu : \pi_2(M, L) \to {\mathbb R}.$$
Let 
$$\Gamma_\omega = \frac{\pi_2(M)}{\ker I_\omega \cap \ker I_c} \text{ and } \Gamma_L = \frac{\pi_2(M, L)}{\ker I_\omega \cap \ker I_\mu},$$
where one
could as well replace $\pi_2(M)$ and $\pi_2(M, L)$ by their images under the respective Hurewicz homomorphisms, namely, the spherical homology groups $H^S_2(M)$ and  $H^S_2(M, L)$, since the quotients are the same\footnote{
If one adopts the point of view in \cite{McDuff} so that \eqref{intro:commutdiag} holds for the non-extended groups, then the $\pi_2$'s are replaced by the respective the spherical homology in $\mathbb R$-coefficients, namely $H^S_2(M;\mathbb R)$ and $H^S_2(M,L;\mathbb R)$. The maps $I_\omega$, $I_c$ and $I_\mu$ are well defined on these homology groups and $\Gamma_\omega$ and $\Gamma_L$ are defined as the respective quotients.
}.
Then the Novikov rings for quantum (or Floer) homology are defined as follows:
$$\Lambda_\omega = \left\{\left.\sum_{B \in \Gamma_\omega} a_B e^{B} \right| a_B \in {\mathbb R} \text{ and } \forall K \in {\mathbb R}, \#\{B | a_B \neq 0 \text{ and } \omega(B) < K\} < \infty \right\}$$
$$\Lambda_L = \left\{\left.\sum_{B \in \Gamma_L} a_B e^{B} \right| a_B \in {\mathbb R} \text{ and } \forall K \in {\mathbb R}, \#\{B | a_B \neq 0 \text{ and } \omega(B) < K\} < \infty\right\}.$$
Their degrees are defined by $\deg(e^{B}) = -2I_c(B)$ and $\deg(e^{B}) = -I_\mu(B)$. 
The homotopy exact sequence $\pi_2(L) \to \pi_2(M) \to \pi_2(M, L) \to \pi_1(L)$ induces an inclusion on the quotients
$$i : \Gamma_\omega \to \Gamma_L,$$
since $\ker I_\omega\cap \ker I_c$ is mapped to $\ker I_\omega \cap \ker I_\mu$. It follows that there is a natural inclusion of the Novikov rings:
$$i : \Lambda_\omega \to \Lambda_L.$$
We can then make a $\Lambda_L$-module into a $\Lambda_\omega$-module via this inclusion.

\subsection{The flow equation}\label{lagFloer:functeq}
Let $H : [0,1] \times M \to {\mathbb R}$ be a time-dependent Hamiltonian function and $\mathbf J = \{J_t\}_{t \in [0,1]}$ a time-dependent $\omega$-compatible almost complex structure. The space of such pairs is $$\mathcal{HJ} = C^\infty([0,1] \times M) \times \mathcal J,$$
where $\mathcal J$ is the space of one-parameter families of $\omega$-compatible almost complex structures.
Let 
$$D^2_+ = \{z\in {\mathbb C} : |z| {\leqslant} 1, \Im z {\geqslant} 0\},$$
$\partial_+$ denote the part of boundary of $D^2_+$ on the unit circle, parametrized by $t \in [0,1]$ as $e^{i\pi t}$, and $\partial_0$ the part on the real line, parametrized by $t \in [0,1]$ as $2t -1$.

We consider the path space
$${\mathcal{P}}_L M = \{l : ([0,1], \{0,1\}) \to (M,L) | [l] = 0 \in \pi_1(M,L)\}.$$
and the covering space $\widetilde {\mathcal{P}}_L M$ of ${\mathcal{P}}_L M$ whose elements are the equivalence classes:
$$[l, w] \text{ where } l \in {\mathcal{P}}_L M \text{ and } w:  (D^2_+; \partial_+, \partial_0) \to (M; l, L),$$
where
$$(l, w) \sim (l', w') \iff l = l' \text{ and } I_\omega(w\#(-w')) = I_\mu(w\#(-w')) = 0.$$
The action functional on $\widetilde {\mathcal{P}}_L M$ is given by
$$a_{H}([l, w]) = -\int_{D^2_+}w^*\omega + \int_{[0,1]} H_t(l(t)) dt,$$
where we use the convention $dH = - \iota_{X_H} \omega$ for the Hamiltonian vector fields. 
An element $\widetilde{l} = [l, w] \in \widetilde {\mathcal{P}}_L M$ is a critical point of $a_H$ if and only if $l$ is a Hamiltonian path connecting points on $L$. 
\begin{definition}\label{lagFloer:nondegen}
A critical point $\widetilde{l}$ is \emph{nondegenerate} if $d\phi_1 T_{l(0)} L \pitchfork T_{l(1)}L$, where $\phi_{t \in[0,1]}$ is the Hamiltonian isotopy generated by $H_{t \in [0,1]}$.
\end{definition}
In a way similar to the case of Hamiltonian Floer homology on $M$, we have
\begin{prop}\label{lagFloer:genric}
For a generic $H$, all critical points of $a_H$ are non-degenerate. \qed
\end{prop}

Floer theory studies the negative gradient flow of $a_H$. Let $(,)_{\mathbf{J}}$ be the metric on ${\mathcal{P}}_L M$ defined by 
$$(\xi, \eta)_{\mathbf{J}} = \int_{[0,1]} \omega(\xi(t), J_t\eta(t)) dt,$$
then the equation of \emph{negative} gradient flow for $a_{H}$ is the following perturbed $\mathbf{J}$-holomorphic equation for $u : {\mathbb R}\times [0,1] \to M$:
\begin{equation}\label{lagFloer:floweq}
\left\{\begin{matrix}
\frac{\partial u}{\partial s} + J_t(u)
 \left(\frac{\partial u}{\partial t} - X_{H_t}(u)\right) = 0 & \text{ for all }
 (s, t) \in {\mathbb R} \times [0,1], \\
u|_{{\mathbb R} \times \{0, 1\}} \subset L
\end{matrix}\right.
\end{equation}
The energy $E(u)$ of a solution $u$ of \eqref{lagFloer:floweq}  with respect to the metric induced by $\mathbf{J}$ is defined as its $s$-energy:   
$$E(u)= \int \left|\frac{\partial u}{\partial s} \right|_t^2ds dt$$
where the $t$-metric on $M$ is ${\langle}\xi,\eta{\rangle}_t = \omega (\xi, J_t \eta)$.
Suppose that all critical points of $a_H$ are non-degenerate, and let $u$ be a finite energy solution, then $l_s(t) := u(s, t)$ converges uniformly to Hamiltonian paths in $C^0$-topology, i.e. $\exists l_{\pm}$ critical points of $a_H$ so that $\lim_{s \to \pm\infty} l_s(t) = l_{\pm}(t)$ uniformly in $t$.

\subsection{Conley-Zehnder index}\label{lagFloer:maslov}
For each nondegenerate critical point $\widetilde{l} = [l, w]$, we can define a Conley-Zehnder index $\mu_{H}(\widetilde{l})$. 
Since $D^2_+$ is contractible, we find a symplectic trivialization $\Phi$ of the bundle $w^*TM$ given by $\Phi_z: T_{w(z)}M \to {\mathbb C}^n$ with standard symplectic structure $\omega_0$ on ${\mathbb C}^n$. We require that $\Phi_r(T_{w(r)}L) = {\mathbb R}^n$ for $r \in  \partial_{0}  D^2_+  \subset D^2_+$, which is possible since $\partial_{0}  D^2_+$ is contractible. Then the linearized Hamiltonian flow $d\phi_t$ along $l$ defines a path of symplectic matrices 
$$E_t = \Phi_{e^{i\pi t}} \circ d\phi_t \circ \Phi^{-1}_1 \in Sp({\mathbb C}^n).$$ 
The Conley-Zehnder index of $\widetilde l$ is defined using the Maslov index of paths of Lagrangian subspaces introduced in Robbin-Salamon \cite{RobbinSalamon1}:
\begin{propdef}\label{lagFloer:intind}
The \emph{Conley-Zehnder index} of $\widetilde l$ is defined as $\mu_H(\widetilde l) = \mu(E_t{\mathbb R}^n, {\mathbb R}^n)$; it satisfies:
\begin{enumerate}
\item $\mu_H(\widetilde l)$ does not depend on the trivialization; 
\item $\mu_H(\widetilde l) + \frac{n}{2} \in {\mathbb{Z}}$;
\item under the deck transformation by $\beta \in \Gamma_L$, we have
$\mu_H(\widetilde l\#\beta) = \mu_H(\widetilde l) + I_\mu(\beta)$.
\end{enumerate}
\end{propdef}

\proof
For $(2)$ see \cite{RobbinSalamon1}, Theorem $2.4$. The rest can be shown similarly as in the case of Hamiltonian loops in $M$.
\qed

\begin{definition}\label{lagFloer:chaingroup}The \emph{Floer chain group} is $FC_*(H) = \oplus_{k} FC_k(H) $ where
$$FC_k(H) := \left\{\left.\sum_{\mu_H(\widetilde l) = k} a_{\widetilde l} \widetilde l \right| a_{\tilde l} \in {\mathbb R} \text{ and } \forall K \in {\mathbb R}, \#\{\widetilde l | a_{\widetilde l} \neq 0, a_H(\widetilde l) < K\} < \infty \right\}.$$
\end{definition}
\noindent
It is easy to see that $FC_*(H)$ is a graded module over the Novikov ring $\Lambda_L$ via:
$$e^B \cdot \widetilde l = \widetilde l \# \beta $$
and we have $$e^B \cdot FC_*(H) \subset FC_{*- \deg(e^B)}(H).$$
We note that by the ring inclusion $i: \Lambda_\omega\to \Lambda_L$, $FC_*(H)$ is also a $\Lambda_\omega$-module.

\subsection{The linearized operator and moduli spaces of flows}\label{lagFloer:linear}
Let us suppose that all critical points of $a_{H}$ are non-degenerate and consider the linearized operator of \eqref{lagFloer:floweq} at a finite energy solution $u$
\begin{equation}\label{linear:op}
D_u \xi = \nabla_{\frac{\partial}{\partial s}} \xi + J_t (u)\nabla_{\frac{\partial}{\partial t}} \xi + \nabla_{\xi}J_t(u) \partial_t u - \nabla_{\xi} \left(J_t (u)X_{H_t}(u)\right),
\end{equation}
where $\xi \in \Gamma(u^*TM; L) = \{\xi \in \Gamma(u^*TM) | \xi|_{{\mathbb R} \times \{0,1\}} \subset TL\}$. Under suitable Banach completion, $D_u: L^p_k(u^*TM; L) \to L^p_{k-1}(u^*TM)$ is Fredholm whose index is the expected dimension of the space of solutions near $u$. 

By \cite{RobbinSalamon2} (Theorem 7.1), the index can be identified as the difference of the Conley-Zehnder indices of the two ends:
\begin{prop}\label{lagFloer:indexdiff}
 Let $\widetilde{{\mathcal{M}}}_{H, \mathbf{J}}(M,L; \widetilde{l}_-, \widetilde{l}_+)$ be the space of all solutions of the equation \eqref{lagFloer:floweq} connecting 
$\widetilde{l}_{-}$ to $\widetilde{l}_{+}$ such that
$[\widetilde l_- \# u \# (- \widetilde l_+)] = 0 \in \Gamma_L$. Its expected dimension is then given by:
$${\rm ind } D_u = \mu_H(\widetilde l_-) - \mu_H(\widetilde l_+).$$
\end{prop}
\qed

The unparametrized moduli space is ${\mathcal{M}}_{H, \mathbf{J}}(M,L; \widetilde{l}_-, \widetilde{l}_+) = \widetilde{{\mathcal{M}}}_{H, \mathbf{J}}(M,L; \widetilde{l}_-, \widetilde{l}_+)/{\mathbb R}$ where the ${\mathbb R}$ action is the shifting of $s$. Thus we have in generic conditions:
$$\dim {\mathcal{M}}_{H, \mathbf J}(M, L; \widetilde l_-, \widetilde l_+) = \mu_H(\widetilde l_-) - \mu_H(\widetilde l_+) - 1.$$

\subsection{Coherent orientations}\label{lagFloer:orient}
We will work over ${\mathbb{Q}}$ or ${\mathbb C}$ instead of ${\mathbb{Z}}_2$. For this reason, we impose the following assumption from now on:
\begin{assumption}\label{lagFloer:relspin}
$L$ is \emph{relatively spin}, i.e. $L$ is orientable and $w_2(L) \in H^2(L; {\mathbb{Z}}_2)$ extends to a class in $H^2(M)$.
\end{assumption}
The above assumption implies that the moduli spaces of holomorphic discs with boundary on $L$ can be canonically oriented with the choice of a relatively spin structure on $L$, i.e.
\begin{itemize}
\item an orientation of $L$,
\item an extension of $w_2(L)$ to $H^2(M)$ and
\item a spin structure on $TL \oplus V |_{L_{(2)}}$, i.e. a trivialization of $TL \oplus V|_{L_{(1)}}$ that extends to $L_{(2)}$,
\end{itemize}
where $L_{(2)}$ is the $2$-skeleton of some triangulation of $L$ and $V$ is an oriented real vector bundle on the $3$-skeleton $M_{(3)}$ of $M$ so that $w_2(V)$ extends $w_2(L)$. It follows that $TL \oplus V |_{L_{(2)}}$ is indeed spin.
Starting from these choices, we may assign to the moduli spaces $\widetilde {\mathcal{M}}_{H, \mathbf{J}}(M,L; \widetilde{l}_-, \widetilde{l}_+)$ a coherent orientation (see for example \cite{FOOO}, \S$44$) in the following way.

First, in order to orient the moduli space of half-tubes $\widetilde {\mathcal{M}}_{H, \mathbf{J}}(M,L; \widetilde{l}_-, \widetilde{l}_+)$, we consider essentially an oriented version of the argument for the PSS \cite{Albers}.
It involves another type of moduli spaces $\widetilde {\mathcal{M}}^\pm_{H^\pm, \mathbf J^\pm}(M, L; \widetilde{l})$ consisting of maps from either the capped strip $Z_-$ or $Z_{+}$ (\cite{Albers}):
\begin{equation}\label{lagFloer:caps}{\mathbb C} \supset Z_\pm = D^2_\mp \cup ({\mathbb R}^\pm \times [0, 1]) \xrightarrow {u^\pm} M \text{ so that } u({\partial Z_\pm}) \subset L.\end{equation}
where this time $D^2_{-}$ denotes the closed left half part of the disk of raduis $1/2$ centered at $1/2 i \in {\mathbb C}$ while $D^2_{+}$ denotes the closed right half part.
The coordinates in $Z_\pm$ is $z = s+it$. Choose and fix $(H^\pm, \mathbf J^\pm)$, a pair of smoothly $z$-dependent Hamiltonian function and almost complex structures so that
$$(H^\pm, \mathbf J^\pm)|_{D^2_\mp} = (0, J) \text{ and } (H^\pm, \mathbf J^\pm)|_{1\mp s < 0} = (H, \mathbf J),$$
where $J$ is a generic almost complex structure on $M$. Consider the equation for $u^\pm : Z_\pm \to M$:
\begin{equation}\label{lagFloer:capeq}
\left\{\begin{matrix}
\frac{\partial u^\pm}{\partial s} + J^\pm_z(u^\pm)
 \left(\frac{\partial u^\pm}{\partial t} - X_{H^\pm_z}(u^\pm)\right) = 0 & \text{ for all }
 (s, t) \in Z_\pm, \\
u|_{\partial Z_\pm} \subset L
\end{matrix}\right.
\end{equation}
The energy $E(u^\pm)$ is defined as the $s$-energy in the usual way, and finite energy solutions $u^\pm$ converge uniformly to Hamiltonian paths of $H$ when $s \to \pm \infty$. Then set
$$\widetilde {\mathcal{M}}_{H^\pm, \mathbf J^\pm}(M, L; \widetilde l) := \left\{u^\pm: Z_\pm \to M \left| \begin{matrix} u^\pm \text{ satisfies \eqref{lagFloer:capeq}},\\ \lim_{s\to \pm\infty} u^\pm = l \text{ and }\\ 
[\widetilde l \# (-u^\pm)] = 0 \in \Gamma_L \end{matrix}\right.\right\}.$$
There are evaluation maps for these moduli spaces, at the points $p_\pm = \pm 1/2 + 1/2 i \in D^2_\pm$:
$$ev^{\pm} : \widetilde {\mathcal{M}}_{H^\pm, \mathbf J^\pm}(M, L; \widetilde l) \to L : u^\pm \mapsto u^\pm(p_\mp).$$

We argue that a choice of the orientations of all the moduli spaces of the form $\widetilde {\mathcal{M}}_{H^+, \mathbf J^+}(M, L; \widetilde l^+)$ induces the orientations of the moduli spaces of the form $\widetilde {\mathcal{M}}_{H^-, \mathbf J^-}(M, L; \widetilde l^-)$ where $l^+ = l^-$. We consider the gluing of the equations
\eqref{lagFloer:capeq} for the moduli spaces $\widetilde {\mathcal{M}}_{H^+, \mathbf J^+}(M, L; \widetilde l^+)$ and $\widetilde {\mathcal{M}}_{H^-, \mathbf J^-}(M, L; \widetilde l^-)$ along $l$. That is, choose and fix an appropriate cut off function $\beta$ and consider the domains
$$Z_{+,R} = D^2_- \cup ([0, R+1] \times [0, 1]) \text{ and } Z_{-,R} = D^2_+ \cup ([-R-1, 0] \times [0,1]),$$
and use $\beta$ to glue the two equations on $Z_\pm$ to define an equation on the glued domain
$$Z_R := Z_{+, R} \sqcup Z_{-, R} / (z \sim z-R-1 \text{ in the ends}).$$
We note that $Z_R$ is conformal to $D^2$ and the equation on $Z_R$ is in fact a compact perturbation of the $\bar\partial_J$-equation for discs with boundary on $L$. 
Because the moduli space of discs is canonically oriented by the choice of a relatively spin structure, we see that the moduli space $\widetilde {\mathcal{M}}_{H^\pm, \mathbf J^\pm}(M, L; \widetilde l, R)$ for the glued equation on $Z_R$ is oriented. From the additivity of indices by standard gluing arguments, we see that the orientations of the $+$-moduli spaces induce orientations of the $-$-moduli spaces. 

Let $B \in \pi_2(M, L)$ and consider $\widetilde l^B = \widetilde l \# B$. When $\widetilde {\mathcal{M}}_{H^+, \mathbf J^+}(M, L; \widetilde l^B)$ is not empty, its orientation is defined from that of $\widetilde {\mathcal{M}}_{H^+, \mathbf J^+}(M, L; \widetilde l)$ and the $\bar\partial$-equation for discs with boundary on $L$ representing class $B$. We note that the moduli space of discs might be empty, or the $\bar \partial$-operator might be non-surjective. Nevertheless, an orientation can be assigned to the index of the $\bar\partial$-operator. Summarizing, we have
\begin{prop}\label{lagFloer:orientcap}
The orientations of the moduli spaces $\widetilde {\mathcal{M}}_{H^\pm, \mathbf J^\pm}(M, L; \widetilde l)$ are determined by the canonical orientations on the indices of the $\bar\partial$-operators of discs with boundaries on $L$ as well as a choice of the orientations on $\widetilde {\mathcal{M}}_{H^+, \mathbf J^+}(M, L; \widetilde l_j)$ for a $\Lambda_L$-basis $\{\widetilde l_j\}$ of $FC_*(H)$.
\qed\end{prop}

\begin{definition}\label{lagFloer:preferredbase}
The basis $\{\widetilde l_j\}$ is called a \emph{preferred basis} for the orientation of the Floer complex $FC_*(M, L; H, \mathbf J)$.
\end{definition}

To obtain the orientations for the moduli spaces $\widetilde {\mathcal{M}}_{H, \mathbf J}(M, L; \widetilde l_-, \widetilde l_+)$, we notice, for example, that  gluing these latter moduli spaces with the moduli spaces $\widetilde {\mathcal{M}}_{H^+, \mathbf J^+}(M, L; \widetilde l_-)$ yields the moduli spaces $\widetilde {\mathcal{M}}_{H^+, \mathbf J^+}(M, L; \widetilde l_+)$. Since both the latter two have been given orientations, these orientations canonically determine orientations on the moduli spaces of half-tubes. Considering the opposite gluing, that is to say using $\widetilde {\mathcal{M}}_{H^-, \mathbf J^-}(M, L; \widetilde l_+)$ instead of $\widetilde {\mathcal{M}}_{H^+, \mathbf J^+}(M, L; \widetilde l_-)$, would give the same induced orientations. 
It is now easy to see that the orientations introduced on $\widetilde {\mathcal{M}}_{H, \mathbf J}(M, L; \widetilde l_-, \widetilde l_+)$ are naturally coherent in the sense of Hofer-Salamon \cite{HoferSalamon}.

\subsection{Floer homology}\label{lagFloer:hlgy} 

From now on, we consider only monotone Lagrangians, i.e. satisfies the following:
\begin{equation}\label{lagFloer:condition}
\text{ there is } \lambda > 0 \text{ such that } I_\omega = \lambda I_\mu \text{ on } \pi_2(M, L).
\end{equation}
Together with assumption \ref{lagFloer:relspin}, we see that the minimal Maslov number of $L$ is at least $2$. The monotonicity condition also ensures that there are no non-trivial  holomorphic spheres with non-positive Chern numbers or non-trivial discs with boundary on $L$ with non-positive Maslov index.

Let $M_k(\mathbf J)$ denote the set of points of $M$ lying on non-constant $J$-holomorphic spheres with Chern number $\leqslant k$, $L_k(\mathbf J)$ the set of points of $L$ lying on the boundary of non-constant $J$-holomorphic discs with Maslov number $\leqslant k$ and $P(H)$ be the set of points of $M$ lying on connecting orbits of $H$.
In the following, we will assume that the pair $(H, \mathbf J)$ is \emph{regular} in the sense that
\begin{itemize}
\item all $J_{0/1}$-holomorphic discs with Maslov index $2$ are regular,
\item $\mathbf J$ is regular for pseudo-holomorphic spheres with Chern number $1$,
\item all connecting orbits of $H$ are non-degenerate,
\item $D_u$ is surjective for finite energy solutions $u$ of \eqref{lagFloer:floweq} with $\text{index } D_u \leqslant 2$.
\item $P(H) \cap M_1(\mathbf J) = \emptyset$ and $P(H) \cap L_2(\mathbf J)$ is empty or of dimension $0$.
\end{itemize}
Standard arguments (cf. e.g. \cite{HoferSalamon}) implies that generic pairs are regular.

The Floer chain complex $FC_*(H, \mathbf{J})$ is 
given by the Floer chain group $FC_*(H)$ with 
the boundary map defined from counting the $0$-dimensional moduli space of solutions:
$$\partial_{H, \mathbf{J}} \widetilde{l}_- = \sum_{\mu_H(\widetilde l_-) = \mu_H(\widetilde l_+) + 1} \#{\mathcal{M}}_{H,\mathbf{J}}(M,L; \widetilde{l}_-, \widetilde{l}_+) \widetilde{l}_+,$$
and extending linearly. 
We then show that
\begin{prop}\label{lagFloer:dsquare}
With assumptions \ref{lagFloer:relspin} 
and assume that $(H, \mathbf J)$ is regular, then $\partial_{H, \mathbf J}^2 = 0$.	
\end{prop}
{\it Proof:}
Writing
$$\partial_{H, \mathbf J}^2 \widetilde l_- = \sum_{\mu_H(\widetilde l_-) = \mu_H(\widetilde l_0) + 1} \#{\mathcal{M}}_{H,\mathbf{J}}(M,L; \widetilde{l}_-, \widetilde{l}_0) \sum_{\mu_H(\widetilde l_0) = \mu_H(\widetilde l_+) + 1} \#{\mathcal{M}}_{H,\mathbf{J}}(M,L; \widetilde{l}_0, \widetilde{l}_+) \widetilde l_+$$
we see that the proposition is equivalent to saying that for each pair of $\widetilde l_-$ and $\widetilde l_+$, we have
$$\sum_{\mu_H(\widetilde l_-) = \mu_H(\widetilde l_0) + 1} \sum_{\mu_H(\widetilde l_0) = \mu_H(\widetilde l_+) + 1} \#{\mathcal{M}}_{H,\mathbf{J}}(M,L; \widetilde{l}_0, \widetilde{l}_+) \#{\mathcal{M}}_{H,\mathbf{J}}(M,L; \widetilde{l}_-, \widetilde{l}_0) = 0.$$
The summand above is the counting for the moduli space of the broken half-tubes connecting $\widetilde l_\pm$. The moduli space of broken half-tubes is part of the boundary components of the $1$-dimensional moduli space ${\mathcal{M}}_{H,\mathbf{J}}(M,L; \widetilde{l}_-, \widetilde{l}_+)$. 

Let $\mathcal C$ be a connected component of the (compactification) of ${\mathcal{M}}_{H,\mathbf{J}}(M,L; \widetilde{l}_-, \widetilde{l}_+)$. A boundary point of $\mathcal C$ is of type $I$ if it is a broken half-tube, is type $II$ if it is a bubbling off of holomorphic discs. The counting in $\partial^2_{H, \mathbf J}$ concerns the type $I$ boundaries. We have the following $3$ cases for $\partial \mathcal C$:
\begin{itemize}
\item empty or is of type $II$ on both ends, or
\item is of type $I$ on both ends, or
\item is of type $I$ on one end and type $II$ on the other.
\end{itemize}
Obviously, if no type $II$ boundary occur in the compactification, an argument similar to the Hamiltonian Floer theory gives the proposition. In the following, we assume that type $II$ boundary does occur. Then the type $I$ boundary and type $II$ boundary are cobordant and the vanishing of counting for either type implies the vanishing of the other. In the following we show the vanishing of counting for the type $II$ boundary points, which would then imply the proposition.

Assume that type $II$ boundary does occur. Then there exist critical points $\widetilde l_\pm$ and a holomorphic disc $v$ with $\mu_L = 2$ so that $v$ is attached to a solution $u$ of \eqref{lagFloer:floweq} such that $\lim_{s \to \pm \infty} u(s, t) = l_\pm$. It follows that $\mu_H(\widetilde l_-) = \mu_H(\widetilde l_+)$. By regularity assumptions, we see that $\widetilde l_- = \widetilde l_+ = \widetilde l$ where $l$ is a connecting Hamiltonian orbit of $H$ and $u(s, t) = l(t)$ for all $s$.
Also by our assumption, $L_2(\mathbf J)$ is compact of dimension $n$. It follows that there are $J_{0/1}$-holomorphic discs through each point of $L$. 

The orientations of the moduli spaces of $J_{0/1}$-holomorphic discs with minimal Maslov number are consistent in the sense that they are connected through cobordisms. On the other hand, 
the orientation of ${\mathcal{M}}_{H, \mathbf J}(M, L; \widetilde l_-, \widetilde l_+) \neq \emptyset$ is obtained as in \S\ref{lagFloer:orient}, by considering the gluing operations, via the canonical orientation of the moduli space of disc together with the choice of orientations on the moduli spaces $\widetilde {\mathcal{M}}_+$. 
The boundary components of ${\mathcal{M}}_{H, \mathbf J}(M, L; \widetilde l_-, \widetilde l_+)$ is oriented by considering the gluing operations. To derive the orientation of the type $II$ boundary points, we need to consider the gluing of the following moduli spaces to the main component $l$:
\begin{quotation}
the moduli spaces $\widetilde {\mathcal{M}}_\pm$ of the capped strips and the moduli space ${\mathcal{M}}_1(M, L; B, J_{0/1})$ of $1$-marked $J_{0/1}$-holomorphic disc with $\mu_L(B) = 2$.
\end{quotation}

The ordering of the gluing operations is given by the orientation of the half-tube ${\mathbb R} \times [0,1]$. Namely, for the case of bubbling off of a disc at $t = 0$, the cyclic order is $\widetilde {\mathcal{M}}_+$, ${\mathcal{M}}_1(M, L; B, J_0)$ then $\widetilde {\mathcal{M}}_-$ and for bubbling off at $t = 1$, the order is $\widetilde {\mathcal{M}}_-$, ${\mathcal{M}}_1(M, L; B, J_1)$ then $\widetilde {\mathcal{M}}_+$.
Note that the orientations of the moduli spaces of holomorphic discs are consistent, while the cyclic ordering of the gluing operations are opposite.
It follows that the counting of configuration of bubbling off at $t = 0$ and $t = 1$ have opposite signs.
It then follows that the counting of type $II$ boundary points vanishes, which implies $\partial_{H, \mathbf J}^2 = 0$.
\qed

Thus 
we define the Floer homology of $(M, L)$ for the regular pair $(H, \mathbf J)$ to be 
$$FH_*(M,L; H, \mathbf{J}) = H_*(FC_*(H, \mathbf{J}), \partial_{H, \mathbf{J}}).$$
The independence of $FH_*(M,L; H, \mathbf{J})$ with respect to the choices of (regular) $H$ 
and $\mathbf{J}$ can be seen using the usual arguments of continuation principle and homotopy 
of homotopies.

\subsection{Half pair of pants product}\label{lagFloer:prod}
\begin{figure}[!h]
\includegraphics[width=0.9\textwidth]{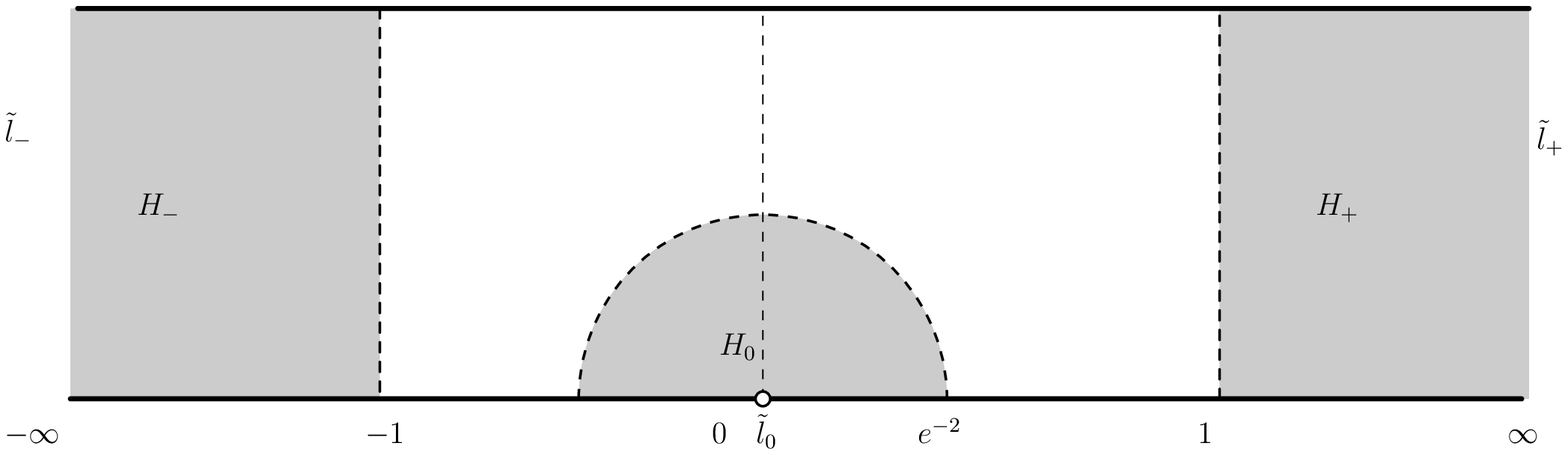}
\caption{}
\label{fig:lagFloer:ppdomain}
\end{figure}

The product on $FH_*(M,L)$ can be defined by ``half pair-of-pants'', perturbed similarly as 
in Seidel \cite{Seidel}, as following. Consider the half cylinder with a boundary puncture $\Sigma_0 = {\mathbb R} \times 
[0,1]\setminus \{(r,0)\}$. The surface $\Sigma_0$ has three ends $e_{\pm}$ and $e_0$:
$$e_+ : [1, \infty) \times [0,1] \to \Sigma_0 \text{ and } e_-: (-\infty, -1] \times [0,1] \to \Sigma_0$$
where $e_{\pm}(s, t) = (s, t)$,
and
$$e_0 : (-\infty, -1] \times [0,1] \to \Sigma_0 : (b, \theta) \mapsto s+it = e^{b - 1+\pi i \theta}$$
is holomorphic with respect to the standard complex structures on the domain and target, 
whose image lies completely in $(-\frac{1}{2}, \frac{1}{2})\times (0, \frac{1}{4})$. 
The ends $e_-$ and $e_0$ are the ``incoming'' ends and $e_+$ is the ``outgoing end''. 
We choose regular pairs $(H_{\pm}, \mathbf{J}_{\pm})$ and $(H_0, \mathbf{J}_0)$ for the 
corresponding ends. Consider the pair $(\mathbf{H}, \mathbf{J})$ where $\mathbf{H} \in C^{\infty}
(\Sigma \times M)$ and $\mathbf{J}$ is a family of compatible almost complex structures 
parametrized by $\Sigma$, such that the pull back of $(\mathbf{H}, \mathbf{J})$ by the maps 
$e_{*}$ is equal to the corresponding pair $(H_*,\mathbf{J}_*)$. Furthermore, we require that
$\mathbf{H}$ restricts to $0$ over $e_0([-2, -1] \times [0,1]) \times M$.
\begin{remark}\label{lagFloer:notion} 
\rm{Here and in the following, a region $\mathscr D \subset {\mathbb R} \times [0,1]$ is provided with cylindrical coordinates
if there is a biholomorphic map $e : I \times S \to \mathscr D$ where $I \subset {\mathbb R}$ is a (possibly infinite) interval
and $S = [0,1]$ or ${\mathbb R}/{\mathbb{Z}}$.
When we ask for the regular pair $(\mathbf{H}, \mathbf J)$ to pull back to
a pair $(H', \mathbf J')$ on a region provided with cylindrical coordinates $I \times S$, we mean that
there is a sequence of (nonempty) smaller intervals 
$$I'' \subset \bar {I''} \subsetneq (I')^\circ \subset \bar{I'} \subsetneq I^\circ \subset I,$$
so that $(\mathbf H, \mathbf J)$ pull-backs to $(H', \mathbf J')$ on $e(I'' \times S)$ while it pulls-back to $(0, J_0)$ on $e((I \setminus I') \times S)$
for some fixed generic compatible almost complex structure $J_0$.
}
\end{remark}
\noindent
The description is conveniently summarized on figure \ref{fig:lagFloer:ppdomain}.
Let $({\mathbb R} \times 
[0,1])^0 = ({\mathbb R} \times [0,1]) \setminus e_0((-\infty, -1] \times [0,1])$ and consider 
the equation 
\begin{equation}\label{lagFloer:pants}
\left\{\begin{matrix}
\frac{\partial u}{\partial s} + J_{s, t}(u)
 \left(\frac{\partial u}{\partial t} - X_{H_{s, t}} (u)\right) = 0 & \text{ for }
 (s, t) \in ({\mathbb R} \times [0,1])^0, \\
\frac{\partial u_0}{\partial s} + J_{e_0(s, t)}(u_0)
 \left(\frac{\partial u_0}{\partial t} - X_{H_{e_0(s, t)}} (u_0)\right) = 0 & \text{ for }
 (s, t) \in (-\infty, -1] \times [0,1], \\
u|_{({\mathbb R} \times \{0, 1\}) \setminus \{(0,0)\}} \subset L
\end{matrix}\right.
\end{equation}
where $u_0 = u \circ e_0$. On the ends $e_*$, a solution $u$ of finite energy again 
limits to critical points $\widetilde l_*$ for the Floer action functional $a_{H_*}$ when 
$s\to \pm\infty$. The half pair-of-pants product is then defined on 
the chain level by counting the $0$-dimensional moduli space ${\mathcal{M}}_{\mathbf H, \mathbf J}
(M, L; \widetilde l_-, \widetilde l_0, \widetilde l_+)$ of such solutions:
$$\widetilde l_- * \widetilde l_0 = \sum_{\widetilde l_+} \#{\mathcal{M}}_{\mathbf H, \mathbf J}(M, L; \widetilde l_-, \widetilde l_0, \widetilde l_+) \widetilde l_+.$$

The orientation of the moduli spaces involved is obtained by considering the gluing with the respective moduli spaces of capped strips. More precisely, gluing with $\widetilde {\mathcal{M}}_{H_-^+,\mathbf J_-^+}(M, L; \widetilde l_-)$ and $\widetilde {\mathcal{M}}_{H_0^+, \mathbf J_0^+}(M, L; \widetilde l_0)$ gives a compact perturbation of the moduli problem for $\widetilde {\mathcal{M}}_{H_+^+, \mathbf J_+^+}(M, L; \widetilde l_{+})$. The order of the gluing operation is first on the end $e_0$ then $e_-$. With all these fixed, we give ${\mathcal{M}}_{\mathbf H, \mathbf J}(M, L; \widetilde l_-, \widetilde l_0, \widetilde l_+)$ the induced orientation. We note that, implicitly, we are also using the orientation of the moduli spaces of holomorphic discs.

To show that it passes to homology, we again look at the boundary of the $1$-dimensional 
moduli space of solutions.
The assumption on the minimal Maslov number implies that a generic family does not have 
any disc bubbling and thus all $1$-dimensional moduli spaces can be compactified by adding broken trajectories.

\begin{remark}\label{lagFloer:productdiffmod}
\rm{
There are two boundary components in ${\mathbb R} \times [0,1]$. In the discussion above, we could also puncture the half cylinder at $(0,1)$ instead and all the arguments will go through and end up with a product of the form $\widetilde l_0 * \widetilde l_-$. The model used above, which is used in this article, will be called the \emph{right} model and the one which punctures $(0,1)$ will be called the \emph{left} model.
}
\end{remark}

\begin{figure}[!h]
\includegraphics[width=0.9\textwidth]{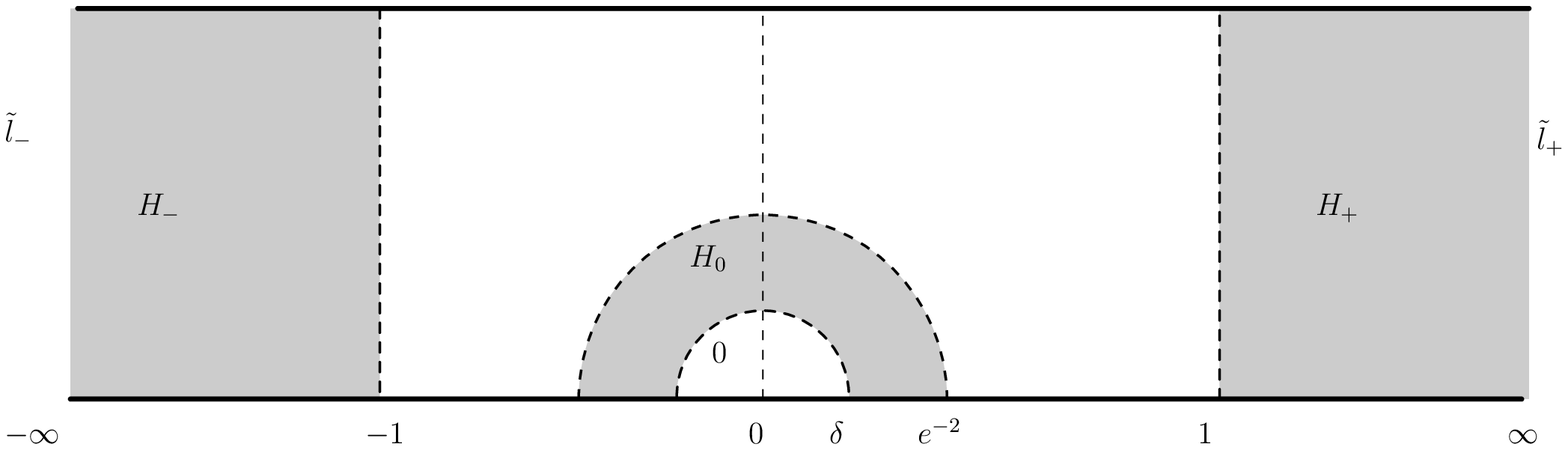}
\caption{}
\label{fig:lagFloer:identitydomain}
\end{figure}

Now we assume that the identity exists for the product just defined and give a description of it in the following. For $\delta \ll 1$, the domain we consider is the unpunctured domain in figure \ref{fig:lagFloer:identitydomain}.
The semi-annulus labelled by $H_0$ is biholomorphic to the half-cylinder $[\ln\delta + 1, -1] \times [0,1]$ by the following:
$$e_\delta : [\ln \delta + 1 , -1] \times [0,1] \to {\mathbb R} \times [0,1] : (b, \theta) \mapsto s+it = e^{b-1 + \pi i \theta}.$$
As $\delta \to 0$, the length of the half-cylinder goes to $\infty$.
We then choose a regular pair $(\mathbf H_\delta, \mathbf J_\delta)$ that pulls-back to the respective regular pairs on the shaded regions and to $(0, J_0)$ on the half disc labelled by $0$. When $\delta > 0$ counting of the $0$-dimensional moduli spaces of the solutions to the perturbed $\bar\partial$ equation described by figure  \ref{fig:lagFloer:identitydomain} gives the connecting homomorphism between $FH_*(M, L; H_-, \mathbf J_-)$ and $FH_*(M, L; H_+, \mathbf J_+)$, which by definition is the identity on the Floer homology $FH_*(M, L)$.

\begin{figure}[!h]
\includegraphics[width=0.6\textwidth]{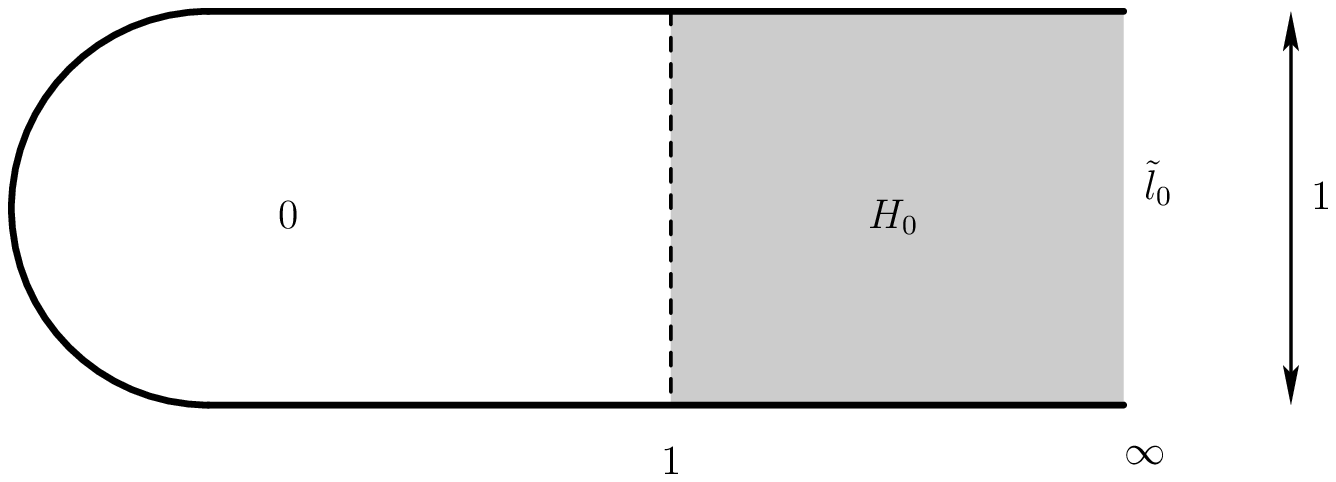}
\caption{}
\label{fig:lagFloer:capdomain}
\end{figure}

Now let $\delta \to 0$, then in the limit the domain splits into two parts, one of which is the domain in figure \ref{fig:lagFloer:ppdomain}. Another part is the ``cap domain'' as shown in figure \ref{fig:lagFloer:capdomain}.
Now let $\widetilde {\mathcal{M}}_{H_0^+, \mathbf J_0^+} (M, L; \widetilde l_0)$ denote the corresponding moduli space of caps (as considered in \S\ref{lagFloer:orient}) defined by the above domain, and consider the element defined from counting dimension $0$ moduli spaces:
$$\idelt_L := \sum_{\widetilde l_0} \#\widetilde {\mathcal{M}}_{H_0^+, \mathbf J_0^+} (M, L; \widetilde l_0) [\widetilde l_0] \in FH_*(M, L).$$
A gluing argument as those in \cite{Albers} then shows that multiplication by $\idelt_L$ gives the identity map of $FH_*(M,L)$, i.e. $\idelt_L$ is the identity of $FH_*(M, L)$ under the half-pair-of-pants product.

\subsection{Action of $FH_*(M)$}\label{lagFloer:quantaction}
Putting in our framework ideas that first appeared in Albers \cite{Albers}, and then in  Biran-Cornea \cite{BiranCornea}, we have the following
\begin{prop}\label{lagFloer:rightmod}
There is a natural action $FH_*(M) \otimes_{\Lambda_\omega} FH_*(M, L) \to FH_*(M, L)$, that exhibits $FH_*(M,L)$ as a right $FH_*(M)$-module, where the Novikov ring $\Lambda_\omega$ is defined over $\mathbb R$.
\end{prop}
We first recall some basic notations from the Hamiltonian Floer homology $FH_*(M)$. We consider the space $\Omega M$ of contractible loops in $M$ and its covering space $\widetilde \Omega M$ which fits into the following covering diagram:
$$\Gamma_\omega \to \widetilde \Omega M \to \Omega M.$$
An element $\widetilde \gamma \in  \widetilde \Omega M$ is an equivalent class of pairs $(\gamma, v)$ where
$$\{\gamma : {\mathbb R}/{\mathbb{Z}} \to M\} \in \Omega M \text{ and } v : (D^2, S^1) \to (M, \gamma)$$
such that $v(e^{2\pi i t}) = \gamma(t)$. The equivalence relation is given by
$$(\gamma, v) \sim (\gamma', v') \text{ whenever } \gamma = \gamma' \text{ and } I_\omega(v\# -v') = I_c(v\# -v') = 0.$$

\begin{figure}[!h]
\includegraphics[width=0.9\textwidth]{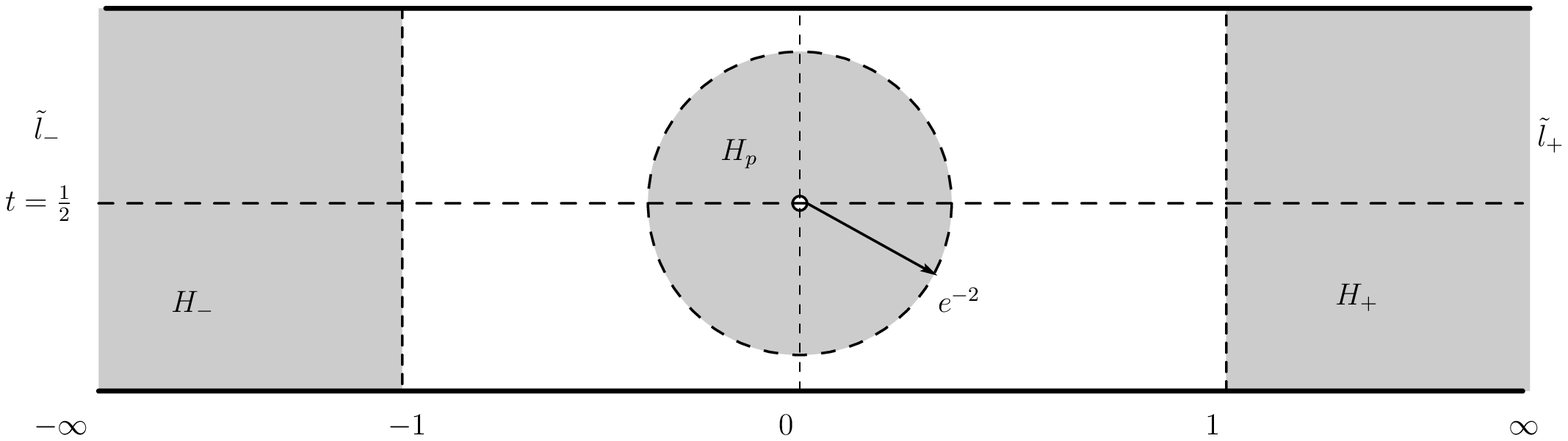}
\caption{}
\label{fig:lagFloer:actiondomain}
\end{figure}

We describe the definition of such an action at the chain level. The domain we consider is $\Sigma_0 = {\mathbb R} \times [0,1] \setminus \{(0, \frac{1}{2})\}$, on which we put the structure of 3 ends, $e_\pm$ as in the previous subsection and $e_p$:
$$e_{p} : (-\infty, -1] \times {\mathbb R}/{\mathbb{Z}} \to \Sigma_0 : (b, \theta) \mapsto s+it = e^{b-1 + 2\pi i \theta} + \frac{i}{2}$$
which is biholomorphic and whose image lies completely in $(-\frac{1}{2}, \frac{1}{2}) \times (\frac{1}{4}, \frac{3}{4})$. We choose regular pairs $(H_\pm, \mathbf J_\pm)$ for $FH_*(M, L)$ on the ends $e_\pm$ and $(H_p, \mathbf J_p)$ for $FH_*(M)$ on the end $e_{p}$. We consider the pair $(\mathbf H, \mathbf J)$ where $\mathbf H \in C^\infty(\Sigma_0 \times M)$ and $\mathbf J$ a family of almost complex structures parametrized by $\Sigma_0$ so that it pull-backs to the respective $(\mathbf H_*, \mathbf J_*)$ on the ends. The domain is again summarized in figure \ref{fig:lagFloer:actiondomain}.
Let $\Sigma_0^p = \Sigma_0 \setminus e_p((-\infty, -1] \times {\mathbb R}/{\mathbb{Z}})$ and consider the equation
\begin{equation}\label{lagFloer:quantacteq}
\left\{\begin{matrix}
\frac{\partial u}{\partial s} + J_{s, t}(u)
 \left(\frac{\partial u}{\partial t} - X_{H_{s, t}} (u)\right) = 0 & \text{ for }
 (s, t) \in \Sigma_0^p, \\
\frac{\partial u_p}{\partial s} + J_{e_p(s, t)}(u_0)
 \left(\frac{\partial u_p}{\partial t} - X_{H_{e_p(s, t)}} (u_p)\right) = 0 & \text{ for }
 (s, t) \in (-\infty, -1] \times {\mathbb R}/{\mathbb{Z}}, \\
u|_{{\mathbb R} \times \{0, 1\}} \subset L
\end{matrix}\right.
\end{equation}
where $u_p = u \circ e_p$. On the ends $e_*$, a solution $u$ of finite energy 
limits to critical points $\widetilde l_{\pm}$ and $\widetilde \gamma_0$ for the respective Floer action functionals when 
$s\to \pm\infty$. The chain level action is then defined 
by counting the $0$-dimensional moduli space ${\mathcal{M}}_{\mathbf H, \mathbf J}
(M, L; \widetilde l_-, \widetilde \gamma_0, \widetilde l_+)$ of such solutions:
$$\widetilde l_- \circ \widetilde \gamma_0 = \sum_{\widetilde l_+} \#{\mathcal{M}}_{\mathbf H, \mathbf J}(M, L; \widetilde l_-, \widetilde \gamma_0, \widetilde l_+) \widetilde l_+.$$
The moduli spaces involved are oriented by the canonical orientation of moduli spaces of discs as well as the choices of orientation for the caps. 
The chain level action passes to homology by the condition $R$, which garantees no bubbling off of holomorphic discs and spheres for the dimension $1$ moduli spaces.
Composing with the PSS isomorphism, we obtain the action of $QH_*(M)$ on $FH_*(M, L)$.

\begin{figure}[!h]
\includegraphics[width=0.9\textwidth]{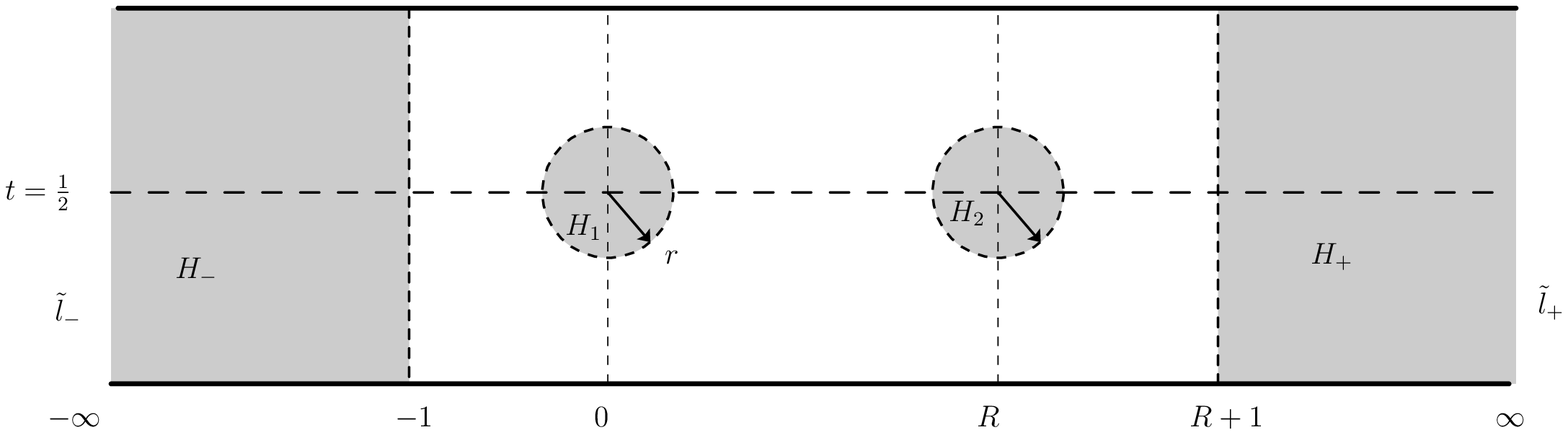}
\caption{}
\label{fig:lagFloer:basicdomain}
\end{figure}

For the sake of completeness and commodity of the reader, we continue to reformulate ideas and results introduced by Albers in our setting: we now show that the action gives a structure of $FH_*(M)$-module. Thus we need here that
\begin{equation}\label{lagFloer:actionmodule}\widetilde l_- \circ (\widetilde \gamma_1 *_{PP} \widetilde \gamma_2) = (\widetilde l_- \circ \widetilde \gamma_1) \circ \widetilde \gamma_2,\end{equation}
where $*_{PP}$ is the pair-of-pants product in $FH_*(M)$.
We consider the twice-punctured domain $\Sigma_{R, 0} = {\mathbb R} \times [0,1] \setminus \{(R, \frac{1}{2}), (0, \frac{1}{2})\}$, where we always set $R > 0$. The basic structure of ends on the domain is illustrated on Figure \ref{fig:lagFloer:basicdomain},
where $r = \epsilon \min (\frac{R}{2}, \frac{1}{2})$ for some $\epsilon \ll 1$. 
The end $e_- : (-\infty, -1] \times [0,1] \to \Sigma_{R,0}$ is the identity map while $e_+ : [1, \infty) \times [0,1] \to \Sigma_{R,0}$ shifts by $R$ to the right.
The structure of the ends in the shaded discs labelled by $H_j$  for $j = 1, 2$ are given by the following
$$e_j : (-\infty, -1] \times {\mathbb R}/{\mathbb{Z}} \to \Sigma_{R,0} : (b, \theta) \mapsto s+it = e^{b+ 2\pi i \theta} + z_j,$$
where $z_1 = \frac{i}{2}$ and $z_2 = R + \frac{i}{2}$, 
and $(H_j, \mathbf J_j)$ are regular pairs for $FH_*(M)$ on the ends $e_j$ for $j = 1, 2$. 
The equation \eqref{lagFloer:actionmodule} is obtained when $R \to 0$ or $R \to \infty$, for which we need compact perturbations of the above basic structure. Because the perturbations are restricted in a compact region of the domain, as well as the condition $R$ excluding bubbling off of discs and spheres, the corresponding operators and resulting moduli spaces (at dimension $0$) allow for a cobordism argument, which will establish the equation \eqref{lagFloer:actionmodule}.

\begin{figure}[!h]
\includegraphics[width=0.9\textwidth]{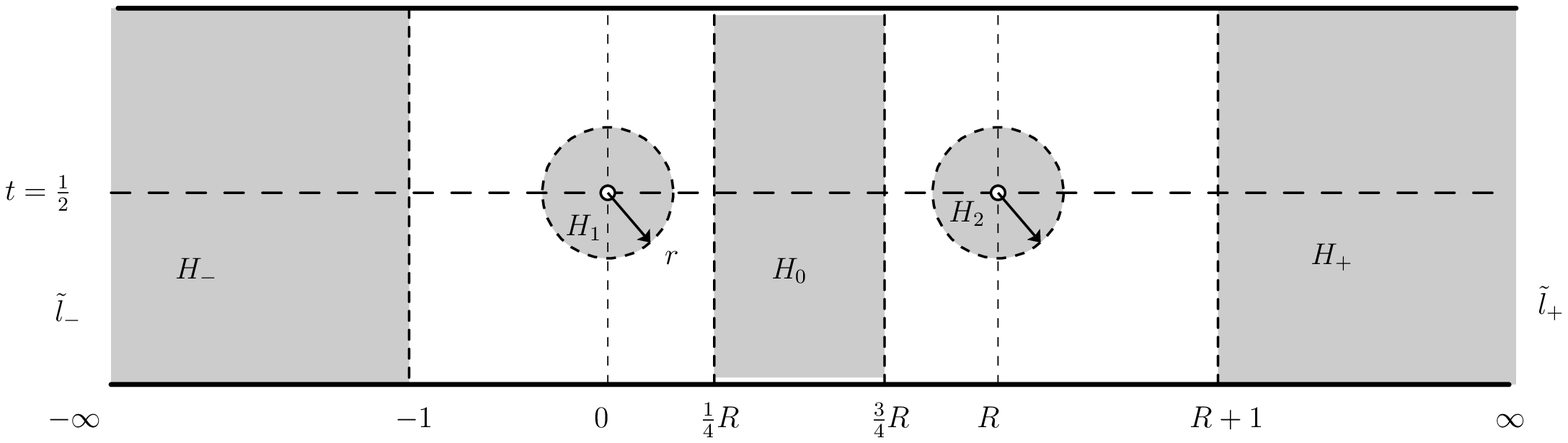}
\caption{}
\label{fig:lagFloer:inftydomain}
\end{figure}

\begin{figure}[!h]
\includegraphics[width=0.9\textwidth]{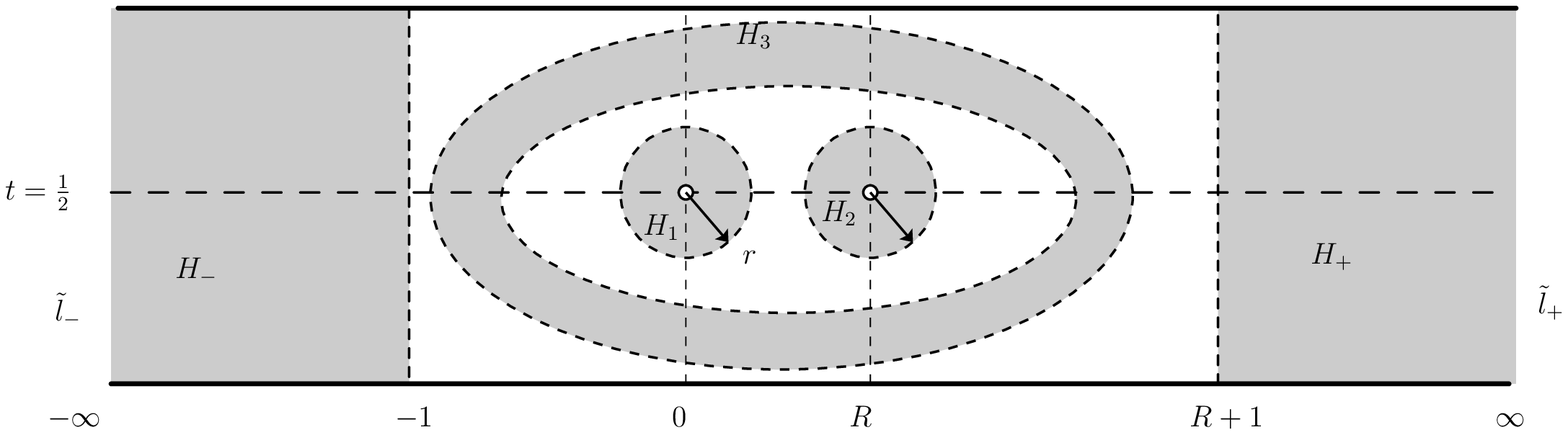}
\caption{}
\label{fig:lagFloer:zerodomain}
\end{figure}

For $R \to \infty$ we choose regular pair $(H_0, \mathbf J_0)$ for $FH_*(M,L)$ and perturb as in figure \ref{fig:lagFloer:inftydomain}.
We write the coordinates explicitly in the region labelled by $H_0$:
$$e_{0,R} : [0, \frac{R}{2}] \times [0,1] \to \Sigma_{R,0} : (s, t) \mapsto (s+\frac{R}{4}, t).$$
When $R \to \infty$, the width of the shaded region labelled by $H_0$ goes to $\infty$. This gives the right hand side of \eqref{lagFloer:actionmodule}.
On the other hand, for $R \ll 1$ we choose a regular pair $(H_3, \mathbf J_3)$ for $FH_*(M)$ and perturb as in figure \ref{fig:lagFloer:zerodomain},
where the shaded region labelled by $H_3$ is an annulus centered at $(\frac{R}{2}, \frac{1}{2})$, for which the outer circle has radius $e^{-1}$ and the inner circle has radius $2R$. This region is biholomorphic to the cylinder $[\ln(2R), -1] \times {\mathbb R}/{\mathbb{Z}}$ via
$$e_{3,R}: [\ln(2R), -1] \times {\mathbb R}/{\mathbb{Z}} \to \Sigma_{R,0} : (b, \theta) \mapsto s+it = e^{b+2\pi i\theta} + \frac{1}{2}(R + i),$$
and we put on the annulus $(H_{3, \theta}, J_{3,\theta})$ using the cylindrical coordinates given by the map $e_{3, R}$. 
When $R \to 0$, the length of the above cylinder goes to $\infty$. This gives the left hand side of \eqref{lagFloer:actionmodule}. 

The equations we consider are similar to \eqref{lagFloer:pants}, where we choose a generic family of pairs $(\mathbf H^R, \mathbf J^R)$ where $\mathbf H^R \in C^\infty(\Sigma_{R, 0} \times M)$ and $\mathbf J^R$ is a family of almost complex structures parametrized by $\Sigma_{R,0}$, so that its pull-back to the cylindrical parts (the shaded parts in the above diagrams) coincides with the corresponding labellings and restricts to $0$ in a neighbourhood of the boundaries of the shaded regions. Then the equation for $(u, R)$ have the form:
\begin{equation}\label{lagFloer:deformeddomains}
\left\{\begin{matrix}
\frac{\partial u}{\partial s} + J^R_{s, t}(u)
 \left(\frac{\partial u}{\partial t} - X_{H^R_{s, t}} (u)\right) = 0 & \text{ for }
 (s, t) \in \text{ unshaded region}, \\
\frac{\partial u_*}{\partial s} + J^R_{e_*(s, t)}(u_*)
 \left(\frac{\partial u_*}{\partial t} - X_{H^R_{e_*(s, t)}} (u_*)\right) = 0 & \text{ for }
 (s, t) \in \text{ the domain of } e_*, \\
u|_{{\mathbb R} \times \{0, 1\}} \subset L
\end{matrix}\right.
\end{equation}
where $u_* = u\circ e_*$ and $* = +, -, 1, 2, (0,R)_{\text{when } R \to \infty}, (3,R)_{\text{when }R \to 0}$ respectively. The gluing that relates the limiting configuration to the configuration where $R \in (0, \infty)$ is similar in all respects to the one employed in the literature, e.g. in \cite{Albers}.

\begin{figure}[!h]
\includegraphics[width=0.5\textwidth]{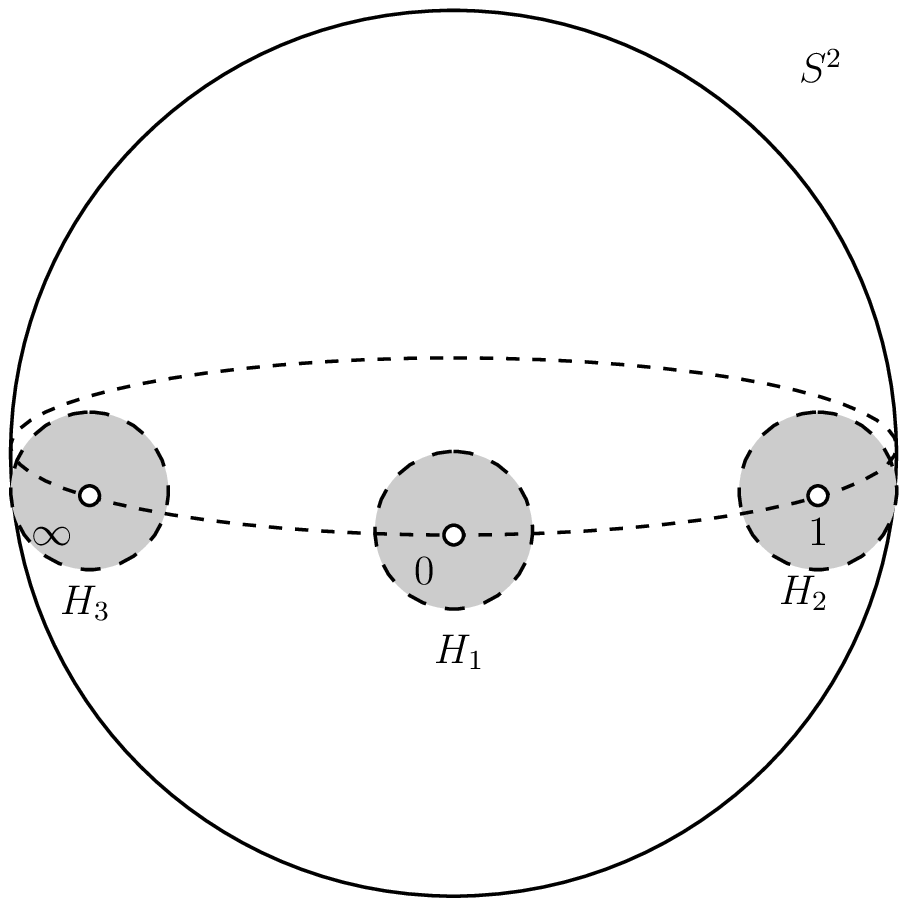}
\caption{}
\label{fig:lagFloer:pantsdomain}
\end{figure}

We recall here that the pair-of-pants product in $FH_*(M, \omega)$ is defined by considering the domain $S^2 \setminus \{0, 1, \infty\}$ (see figure \ref{fig:lagFloer:pantsdomain}).
The shaded discs around $0$, $1$ and $\infty$ are provided with cylindrical coordinates, so that the disc around $0$ and $1$ are identified with $(-\infty, -1] \times {\mathbb R}/{\mathbb{Z}}$ and the one around $\infty$ is identified with $[1, \infty) \times {\mathbb R}/{\mathbb{Z}}$. To set the order of the multiplication, let $\widetilde \gamma_j$ be critical points for the action functional of $H_j$  for $j = 1, 2$. Then counting $0$-dimensional moduli spaces of maps from the above domain satisfying the perturbed holomorphic equation as described in figure \ref{fig:lagFloer:pantsdomain} defines the product $\widetilde \gamma_1 * \widetilde \gamma_2$. Comparing this with the description of figure \ref{fig:lagFloer:zerodomain}, we see that the $FH_*(M)$-action is indeed a right action.

\begin{remark}\label{lagFloer:quantactionrmk}
\rm{Via the PSS-isomorphism between $FH_*(M)$ and $QH_*(M)$, the action described in this section can also be thought of as a right action of $QH_*(M)$ on $FH_*(M, L)$. In a way similar to the description of the PSS-isomorphism, the action by $QH_*(M)$ can be constructed directly using Morse trajectories as in Biran-Cornea \cite{BiranCornea}. 
}
\end{remark}

\subsection{Albers' map}\label{lagFloer:Albers}
We describe a proof of the following
\begin{prop}\label{lagFloer:compatAlbers}
The action introduced above is compatible with the comparison map $\mathscr A : FH_*(M) \to FH_*(M,L)$ introduced by Albers in \cite{Albers}
via the half-pair-of-pants product, whenever all ingredients are defined. It means
$$[\widetilde l_-] \circ [\widetilde \gamma_0] = [\widetilde l_-] * \mathscr A([\widetilde \gamma_0]),$$
where $[\widetilde l_-] \in FH_*(M, L)$ and $[\widetilde \gamma_0] \in FH_*(M)$.
\end{prop}

\begin{figure}[!h]
\includegraphics[width=0.9\textwidth]{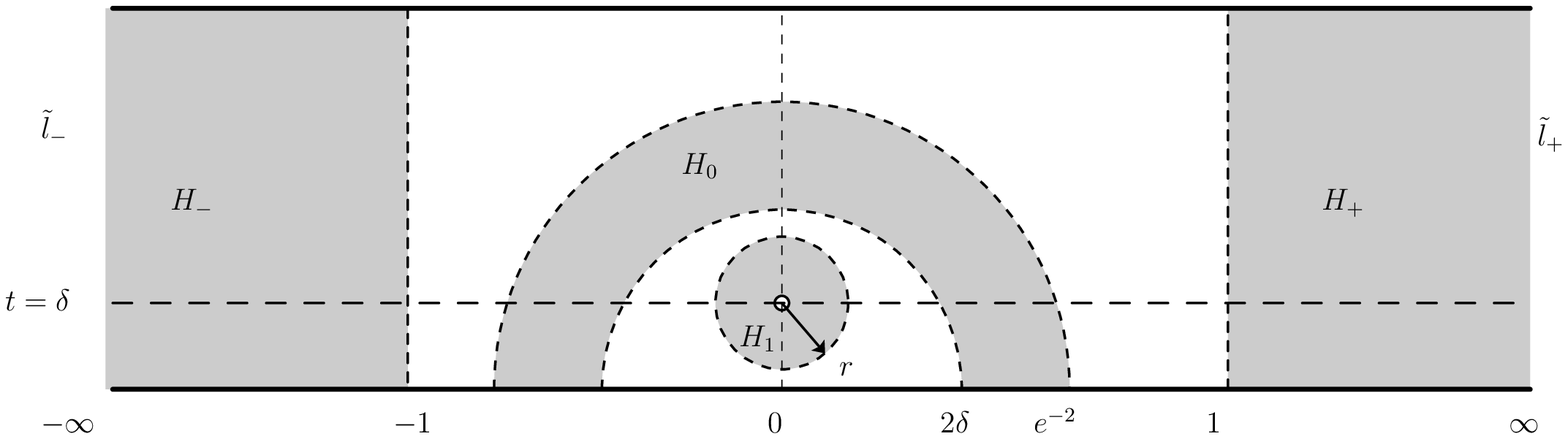}
\caption{}
\label{fig:lagFloer:albersdomain}
\end{figure}

\begin{figure}[!h]
\includegraphics[width=0.6\textwidth]{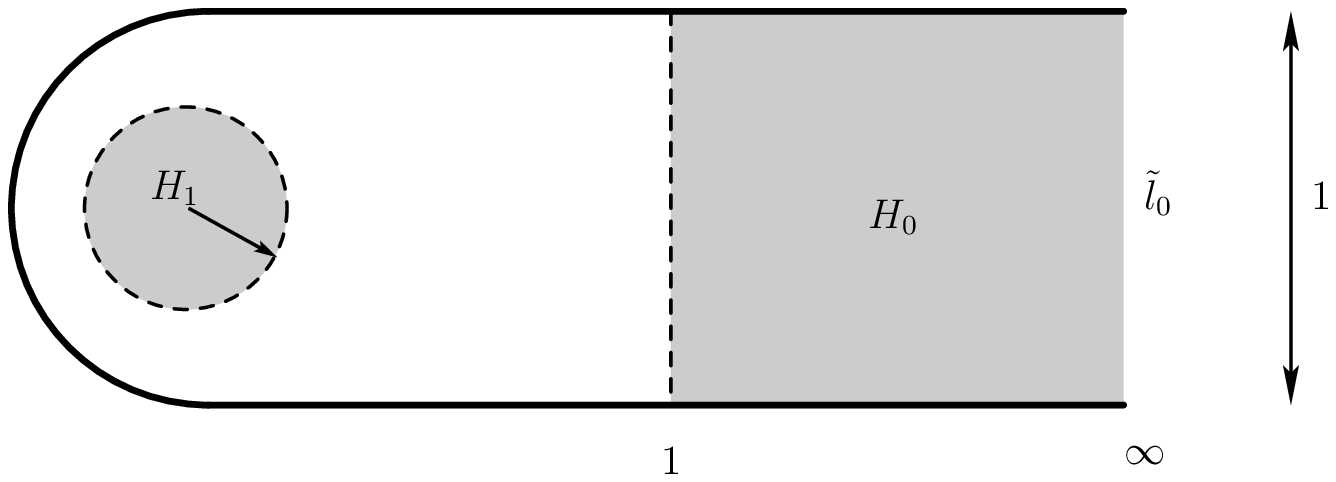}
\caption{}
\label{fig:lagFloer:chimneydomain}
\end{figure}

We consider the domain $\Sigma_{\delta} = {\mathbb R} \times [0,1] \setminus \{(0, \delta)\}$, $\delta \to 0$ with the cylindrical structure as in figure \ref{fig:lagFloer:albersdomain}, where $r = \epsilon \delta$ for some $\epsilon \ll 1$.
We choose a regular pair $(H_0, \mathbf J_0)$ for $FH_*(M, L)$ and $(H_1, \mathbf J_1)$ for $FH_*(M)$. The shaded half-annulus labelled by $H_0$ is centered at $(0,0)$ and outer radius  $e^{-2}$ and inner radius is $2\delta$. The cylindrical coordinates on the shaded half-annulus is given by the biholomorphic map
$$e_{a, \delta} : [\ln (2\delta) + 1, -1] \times [0,1] \to \Sigma_{\delta} : (b, \theta) \mapsto s+it = e^{b -1 + \pi i \theta},$$
and we put on it $(H_{a,\theta}, J_{a, \theta})$ using $e_{a, \delta}$. When $\delta \to 0$, the length of the cylinder goes to $\infty$.
We then solve an equation for $(u, \delta)$ of the type \eqref{lagFloer:deformeddomains} with the above domain and the cylindrical data. In the limit $\delta \to 0$, the domain splits into the domain for the half-pair-of-pants in \S\ref{lagFloer:prod} together with the ``chimney domain'' used to define the map $\mathscr A$, see figure \ref{fig:lagFloer:chimneydomain}.
Then a gluing argument similar to the one in \cite{Albers} proves the statement.

\begin{figure}[!h]
\includegraphics[width=0.9\textwidth]{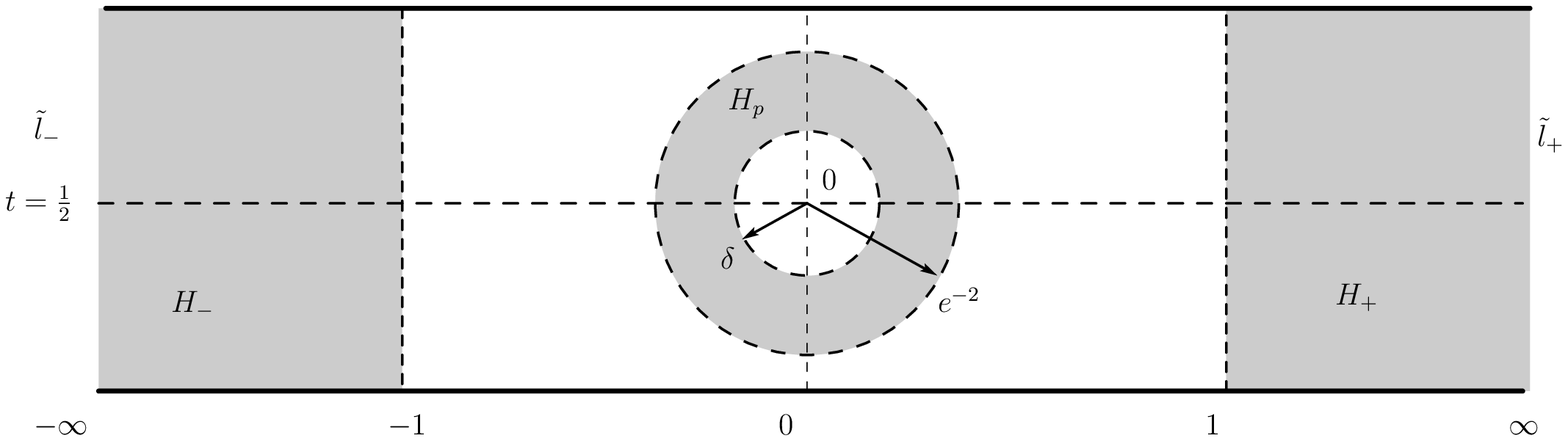}
\caption{}
\label{fig:lagFloer:pssdomain}
\end{figure}

\begin{figure}[!h]
\includegraphics[width=0.5\textwidth]{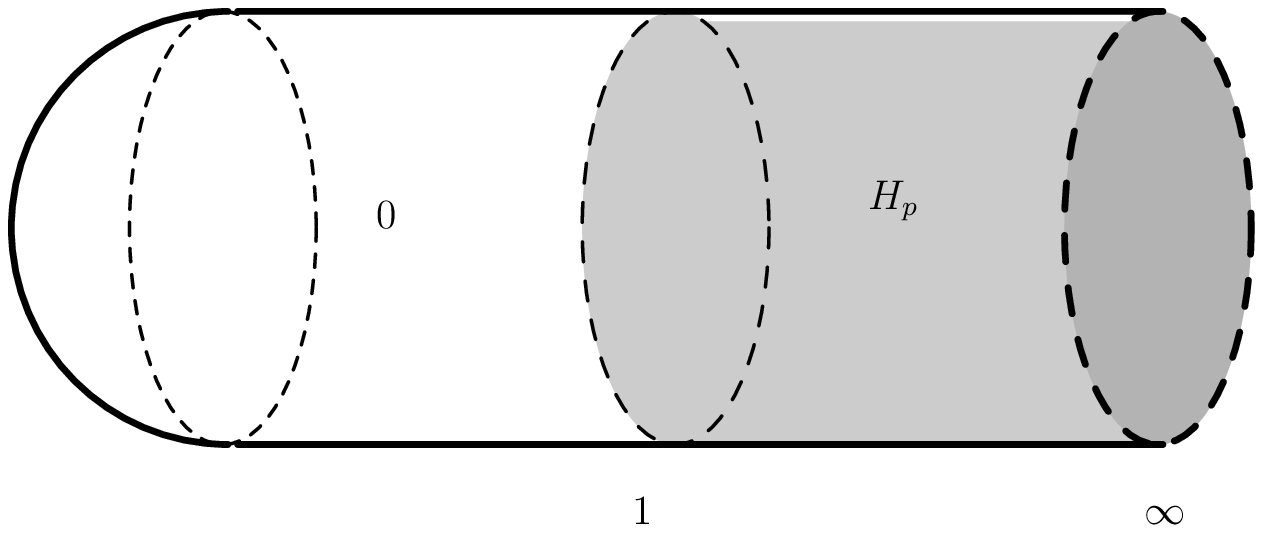}
\caption{}
\label{fig:lagFloer:psscapdomain}
\end{figure}

The following corollary of proposition \ref{lagFloer:compatAlbers} gives the image of the identity $\idelt \in FH_*(M)$ under the comparison map $\mathscr A$:
\begin{prop}\label{lagFloer:albersidentity}
Suppose that identity exists for the half-pair-of-pants product defined in \S\ref{lagFloer:prod}, then $\mathscr A(\idelt) = \idelt_L$.
\end{prop}
{\it Proof:}
According to proposition \ref{lagFloer:compatAlbers}, we only have to show that the action of $\idelt$ on $FH_*(M, L)$ gives the identity map. We consider the domain in figure \ref{fig:lagFloer:pssdomain}, for $\delta \ll 1$.
When $\delta > 0$, this again leads to the identity map on the Lagrangian Floer homology $FH_*(M, L)$. In the limit when $\delta \to 0$,  the domain splits into two parts, one of which is described in figure \ref{fig:lagFloer:actiondomain}. The other one is the ``capped domain'' in the description of PSS map in \cite{PSS}, see figure \ref{fig:lagFloer:psscapdomain}.
A similar argument as in \S\ref{lagFloer:prod} shows that counting the dimension $0$ moduli spaces of caps gives $\idelt \in FH_*(M)$, and a gluing argument then shows that the action of $\idelt$ on $FH_*(M, L)$ is indeed the identity map.
\qed
\section{Seidel's construction for Lagrangian Floer homology} \label{lagSeidel}

For the discussion in this section, we impose furthermore the following
\begin{assumption}\label{lagSeidel:assume}
$FH_*(M,L)$ is non-vanishing, in particular, $L$ is non-displaceable in $M$.
\end{assumption}
\noindent
By this assumption, any (time-dependent) Hamiltonian function has contractible flow line connecting points on $L$.

\subsection{Path group and action}\label{lagSeidel:groupact}
Let ${\rm Ham}_L(M, \omega)$ be the subgroup of ${\rm Ham}(M, \omega)$ that preserves $L$, i.e.
$$\phi \in {\rm Ham}_L(M, \omega) \iff \phi \in {\rm Ham}(M, \omega) \text{ and } \phi(L) = L.$$
We consider the following path group in ${\rm Ham}(M, \omega)$:
\begin{definition}\label{lagSeidel:definition}
${\mathcal{P}}_L{\rm Ham}(M, \omega) := \{g: ([0,1]; \{0\}, \{1\}) \to ({\rm Ham}(M, \omega); id, {\rm Ham}_L(M, \omega))\}$ is a group with pointwise composition:
$$(g \circ h)_t = g_t h_t.$$
For $g \in {\mathcal{P}}_L{\rm Ham}(M,\omega)$, the action of it on a path $l : ([0,1]; \{0,1\}) \to (M, L)$ is
$$(g\circ l)(t) = l^g(t) := g_t\circ l(t).$$
\end{definition}

Suppose that $g$ is generated by $K : [0,1] \times M \to {\mathbb R}$ then we have $g^*\alpha_{H} = \alpha_{H^g}$ and $g^*(,)_{\mathbf J} = (,)_{\mathbf J^g}$ where
$$H^g(t,x) = H(t, g_tx) - K(t,g_tx) \text{ and } J_t^g(x) = dg_t^{-1} \circ J_t(g_tx) \circ dg_t.$$
In particular, we have $(H^g)^h = H^{gh}$ and $(\mathbf J^g)^h = \mathbf J^{gh}$. 
Let $\phi_t$ denote the Hamiltonian isotopy generated by $H_t$, then $H^g_t$ generates $g^{-1}_t \phi_t$. It follows that the connecting Hamiltonian flow lines of $H_t$ and $H_t^g$ are in one-to-one correspondence.

   As in Lalonde-McDuff-Polterovich \cite{LalondeMcDuffPolterovich} where it is shown that ${\rm Ham}(M)$ acts trivially on homology (and sends contractible loops in $M$ to contractible loops), one sees easily that the same argument shows that the action of ${\mathcal{P}}_L{\rm Ham}(M, \omega)$ on the space of paths preserves the component ${\mathcal{P}}_L M$.

Most computations in the following are parallel to the corresponding ones in \cite{Seidel}.
\begin{prop}\label{lagSeidel:lifting}
The action of ${\mathcal{P}}_L{\rm Ham}(M,\omega)$ on ${\mathcal{P}}_L M$ can be lifted to an action of $\widetilde{{\mathcal{P}}}_L{\rm Ham}(M,\omega)$ on the covering $\widetilde{{\mathcal{P}}}_L M$, where 
$$\widetilde{{\mathcal{P}}}_L{\rm Ham}(M,\omega) := \left\{\left.(g, \widetilde g) \in {\mathcal{P}}_L{\rm Ham}(M,\omega) \times {\rm Homeo}(\widetilde{{\mathcal{P}}}_L M) \right| \widetilde{g} \text{ lifts the action of } g \right\}$$
\end{prop}

\proof 
We only need to show that the action can be lifted. Suppose $\gamma : S^1 \to  {\mathcal{P}}_LM$ is a loop that can be lifted to $\widetilde {\mathcal{P}}_LM$, then it is represented by a map $B: \left(S^1 \times [0,1], S^1 \times \{0,1\}\right) \to (M, L)$, such that $\omega(B) = \mu_L(B) = 0$. The loop $\gamma^g = \{g(\gamma_s)\}_{s \in S^1}$ is represented by $B^g(s, t) = g_t \circ B(s, t)$. Because $dg_t : (B^*TM, \partial B^*TL) \to ((B^g)^*TM, (\partial B^g)^*TL)$ is an isomorphism of symplectic bundles preserving the Lagrangian boundary conditions, it follows that $\mu_L(B^g) = \mu_L(B) = 0$.
We compute $(B^g)^*\omega = \omega\left(\frac{\partial B^g}{\partial s}, \frac{\partial B^g}{\partial t}\right) ds\wedge dt = B^*\omega + d\theta$ with $\theta = K(t, B^g(s, t)) dt$. Since $\theta|_{\partial(S^1 \times [0,1])} = 0$ we find that $\omega(B^g) = \omega(B) = 0$. Thus, $\gamma^g$ again can be lifted to $\widetilde {\mathcal{P}}_LM$, which implies that the action of $g$ can be lifted.
\qed

The groups fit into the exact sequence:
$$0\to \Gamma_L \to \widetilde{{\mathcal{P}}}_L{\rm Ham}(M,\omega) \to {\mathcal{P}}_L{\rm Ham}(M,\omega) \to 0,$$
and passing to homotopy, we get the exact sequence:
$$\Gamma_L \to \widetilde \pi_1({\rm Ham}(M, \omega), {\rm Ham}_L(M, \omega)) \to \pi_1({\rm Ham}(M, \omega), {\rm Ham}_L(M, \omega)) \to 0.$$
Let $\widetilde l = [l, w] \in \widetilde {\mathcal{P}}_L M$ (not necessarily a critical point of any functional) and $\Phi_z : T_{w(z)} M \to {\mathbb C}^n$ be any trivialization satisfying $\Phi_r : T_{w(r)}L \to {\mathbb R}^n$ for $r \in [-1, 1]$. Let $\Phi^{\widetilde g}_z$ be a similar trivialization defined for $\widetilde l^{\widetilde g} = [l^g, w^{\widetilde g}]$.
Now consider 
$$G_t = \Phi^{\widetilde g}_{e^{i\pi t}} \circ dg_t \circ \Phi_{e^{i\pi t}}^{-1} : {\mathbb C}^n \to {\mathbb C}^n.$$
Then $\{G_t{\mathbb R}^n\}$ is a loop of Lagrangian subspaces in ${\mathbb C}^n$. The following definition then does not depend on either the trivialization or the choice of $\widetilde l$.
\begin{definition}\label{lagSeidel:actionind}
The \emph{Maslov degree} of $\widetilde g$ is $\mu(\widetilde g) = \mu(G_t{\mathbb R}^n)$. \qed
\end{definition}
\begin{prop}\label{lagSeidel:index}
Let $\widetilde{l} = [l, w]$ be a critical point of $a_{H^g}$, then $\widetilde{l}^{\widetilde{g}}$ is a critical point of $a_H$. Furthermore, $\widetilde l$ is non-degenerate if and only if $\widetilde l^{\widetilde g}$ is so. For such critical points, we have $\mu(\widetilde g) =  \mu_H(\widetilde{l}^{\widetilde{g}}) - \mu_{H^g}(\widetilde{l})$. It follows that $\mu : \widetilde g \mapsto \mu(\widetilde g)$ defines a group homomoprhism $\mu: \widetilde \pi_1({\rm Ham}(M, \omega), {\rm Ham}_L(M, \omega)) \to {\mathbb{Z}}$.
\end{prop}
\proof 
A direct computation from the definitions establish the first statement of the proposition. 
Suppose $l$ is non-degenerate, then $dg_1^{-1} \circ d\phi_1 (T_{l(0)}L) \pitchfork T_{l(1)}L \iff d\phi_1 (T_{l(0)}L) \pitchfork T_{l^g(1)}L$ since $dg_1$ preserves $TL$, and so the second statement follows.

Let $\widetilde l^{\widetilde g} = [l^g, w^g]$ and $\Phi_z : T_{w^g(z)}M \to {\mathbb C}^n$ be a trivialization that defines $\mu_H(\widetilde l^{\widetilde g})$, then
$$\mu_H(\widetilde l^{\widetilde g}) = \mu(E_t{\mathbb R}^n, {\mathbb R}^n), \text{ where } E_t = \Phi_{e^{i\pi t}}\circ d\phi_t \circ \Phi_1^{-1} : {\mathbb C}^n \to {\mathbb C}^n.$$
Let $\Phi^g_z : T_{w(z)}M \to {\mathbb C}^n$ be the trivialization defining $\mu_{H^g}(\widetilde l)$, then 
$$\mu_{H^g}(\widetilde l) = \mu(E^g_t{\mathbb R}^n, {\mathbb R}^n), \text{ where } E^g_t = \Phi^g_{e^{i\pi t}}\circ dg^{-1}_t \circ d\phi_t \circ \left(\Phi^g_1\right)^{-1} : {\mathbb C}^n \to {\mathbb C}^n.$$
Suppose $\Phi^g_1 = \Phi_1$ and let $G_t^{-1} = \Phi^g_{e^{i\pi t}}\circ dg^{-1}_t \circ \Phi_{e^{i\pi t}}^{-1}$, then $E_t = G_t \circ E^g_t = G_t \# E^g_t$ because $G_t$ is a loop. Thus the property of the Maslov index of Lagrangian paths (see \cite{RobbinSalamon1} theorem $2.3$) gives
$$\mu(E_t{\mathbb R}^n, {\mathbb R}^n) = \mu(G_t) + \mu(E^g_t {\mathbb R}^n, {\mathbb R}^n) \Rightarrow \mu(\widetilde g) =  \mu_H(\widetilde{l}^{\widetilde{g}}) - \mu_{H^g}(\widetilde{l}).$$
\qed

In a way entirely parallel to \cite{Seidel}, we have
\begin{prop}\label{lagSeidel:fhaction}
For critical points $\widetilde l_-$, $\widetilde l_+$ of $a_{H^g}$, there is a bijection of moduli spaces:
$$\begin{matrix}
{\mathcal{M}}_{H^g, \mathbf J^g}(M, L; \widetilde l_-, \widetilde l_+) & \to & {\mathcal{M}}_{H, \mathbf J}(M, L; \widetilde l_-^g, \widetilde l_+^g)\\
u & \mapsto & u^g
  \end{matrix}
$$
where 
\begin{equation}\label{lagSeidel:flowaction}
u^g(s, t) := g_t \circ u(s, t).
\end{equation}
Furthermore, $(H, \mathbf J)$ is regular iff $(H^g, \mathbf J^g)$ is.
The map $FC_*(\widetilde g; H, \mathbf J)$ defined by ${\langle}\widetilde l{\rangle} \mapsto {\langle}\widetilde l^g{\rangle}$ passes to homology:
$$FH_*(\widetilde g) : FH_*(H^g, \mathbf J^g) \to FH_{* + \mu(\widetilde g)}(H, \mathbf J)$$
and defines an automorphism of $FH_*(M, L)$ of degree $\mu(\widetilde g)$. Furthermore the following hold:
\begin{enumerate}
\item for $(g, \widetilde g) = (id, id)$, $FH_*(\widetilde g) = id$,
\item for $(g, \widetilde g) = (id, \beta)$ with $\beta \in \Gamma_L$, we have $FH_*(\widetilde g) = \beta \cdot id$,
\item $FH_*(\widetilde g)$ is a $\Lambda_L$-module automorphism of degree $\mu(\widetilde g)$,
\item $FH_*(\widetilde g \circ \widetilde g') = FH_*(\widetilde g)\circ FH_*(\widetilde g')$.
\end{enumerate}
\qed
\end{prop}

\subsection{Homotopy invariance}\label{lagSeidel:htpyinv}
We consider a smooth path $\{g_r\}_{r \in [0,1]}$ starting from the identity in ${\mathcal{P}}_L{\rm Ham}(M, \omega)$ and a lift $\{(g_r, \widetilde g_r)\}$ of it to $\widetilde {\mathcal{P}}_L{\rm Ham}(M, \omega)$. Then proposition \ref{lagSeidel:index} implies that $\mu(\widetilde g_r) = 0$ for all $r \in [0,1]$.
The path $\{g_r\}$ corresponds to a smooth family $\{g_{r,t}\}_{(r,t) \in [0,1]^2}$ of Hamiltonian diffeomorphisms in ${\rm Ham}(M, \omega)$ so that
$$g_{0,t} = g_{r, 0} = id \text{ and } g_{r, 1} \in {\rm Ham}_L(M, \omega) \text{ for all } r \in [0,1].$$
Choose a smooth family of Hamiltonians $K_r : [0,1] \times M \to {\mathbb R}$ for $r \in [0,1]$ so that $K_r$ generates $g_r$ and $K_0 = 0$. Let
$$H^{g_r}(t, x) = H(t, g_{r,t}(x)) - K_r(t, g_{t,r}(x)) \text{ and } J_t^{g_r} = dg_{r, t}^{-1} \circ J_t(g_{r, t}(x)) \circ dg_{r, t}$$

Let $(H_0, \mathbf J_0) = \{(H_t, J_t)\}_{t \in [0,1]}$ and $(H_1, \mathbf J_1) = \{(H^{g_1}_t, J^{g_1}_t)\}_{t \in [0,1]}$, then the construction in the last subsection gives
$$FH_*(\widetilde g_1) : FH_*(H_1, \mathbf J_1) \to FH_*(H_0, \mathbf J_0).$$
Let $(\bar H, \bar{\mathbf J}) = \{(H_{s,t}, J_{s,t})\}_{(s,t) \in {\mathbb R} \times [0,1]}$ be a regular homotopy connecting $(H_0, \mathbf J_0)$ and $(H_1, \mathbf J_1)$:
$$(\bar H_s, \bar{\mathbf J}_s) = \left\{\begin{matrix}(H_1, \mathbf J_1) & s {\leqslant} -1 \\ (H_0, \mathbf J_0) & s {\geqslant} 1 \end{matrix}\right. .$$
We consider the moduli spaces of the solutions of the following equation for maps $u : {\mathbb R} \times [0,1] \to M$ with $\partial u : {\mathbb R} \times \{0,1\} \to L$:
$$\frac{\partial u}{\partial s} + J_{s, t}(u)
 \left(\frac{\partial u}{\partial t} - X_{H_{s, t}} (u)\right) = 0.$$
Here $(\bar H, \bar{\mathbf J})$ being regular means that all solutions $u$ are regular, i.e. their linearizations are surjective.
The moduli space ${\mathcal{M}}_{\bar H, \bar{\mathbf J}}(M, L; \widetilde l_-, \widetilde l_+)$ denotes the space of solutions $u$ that converge to Hamiltonian paths when $s \to \pm \infty$:
$$\lim_{s\to -\infty} = \widetilde l_-^{g_1^{-1}} \text{ and } \lim_{s\to +\infty} = \widetilde l_+,$$
where $\widetilde l_\pm$ are critical points of $a_{H}$. The dimension of the moduli space is given by
$$\mu_H(\widetilde l_-) - \mu_H(\widetilde l_+),$$
since there is no ${\mathbb R}$-action anymore.
Then the continuation map on the chain level
$$\Phi_{\bar H, \bar{\mathbf J}}: FC_*(H_1) \to FC_*(H_0)$$
is defined by counting dimension $0$ moduli spaces:
$$\Phi_{\bar H, \bar{\mathbf J}}(\widetilde l_-^{g_1^{-1}}) = \sum_{\widetilde l_+} \#{\mathcal{M}}_{\bar H, \bar{\mathbf J}}(M, L; \widetilde l_-, \widetilde l_+) \widetilde l_+.$$
That $\Phi_{\bar H, \bar{\mathbf J}}$ is a chain map is shown by considering the dimension $1$ moduli spaces. Thus we have the continuation map for Floer homology, which is also denoted $\Phi_{\bar H, \bar{\mathbf J}}$. The homotopy invariance of $FH_*(\widetilde g)$ is equivalent to the following:

\begin{prop}\label{lagSeidel:deformhtpy}
With the above setup, we have $FH_*(\widetilde g_1) = \Phi_{\bar H, \bar{\mathbf J}}$.
\end{prop}
As in \cite{Seidel}, we consider the deformation of homotopies, from the trivial homotopy to $(\bar H, \bar{\mathbf J})$ by the curve $\{(g_r, \widetilde g_r)\}$, which is a family $(\widetilde H, \widetilde{J}) = \{(H_{r,s,t}, J_{r,s,t})\}_{(r,s,t) \in [0,1] \times {\mathbb R} \times [0,1]}$ where
\begin{eqnarray*}
H_{r,s,t}(x) = H^{g_r}_t(x), & J_{r,s,t}(x) = J^{g_r}_t(x) & \text{ for } s {\leqslant} -1\\
H_{r,s,t}(x) = H_t(x), & J_{r,s,t}(x) = J_t(x) & \text{ for } s {\geqslant} 1\\
H_{0,s,t}(x) = H_t(x), & J_{0,s,t}(x) = J_t(x) & \text{ and } \\
H_{1,s,t}(x) =  H_{s,t}(x) & J_{1,s,t}(x) = J_{s,t}(x) & 
\end{eqnarray*}
The equation that we are now concerned with is the following:
\begin{equation}\label{lagSeidel:homotopyeq}\frac{\partial u}{\partial s} + J_{r, s, t}(u)
 \left(\frac{\partial u}{\partial t} - X_{H_{r, s, t}} (u)\right) = 0,
\end{equation}
for the pair $(r, u)$, where $r \in [0,1]$ and $u: {\mathbb R} \times [0,1] \to M$ with $\partial u : {\mathbb R} \times \{0,1\} \to L$. Let ${\mathcal{M}}_{\widetilde H, \widetilde{\mathbf J}}(M, L; \widetilde l_-, \widetilde l_+)$ denote the moduli space of solutions $(r, u)$ so that $u$ solves the equation at the parameter $r$ and converges to Hamiltonian paths as $s \to \pm \infty$, i.e.:
$$\lim_{s\to -\infty} = \widetilde l_-^{g_r^{-1}} \text{ and } \lim_{s\to +\infty} = \widetilde l_+,$$
where $\widetilde l_\pm$ are critical points of $a_H$. Then the expected dimension of this moduli space is
$$\mu_H(\widetilde l_-) - \mu_H(\widetilde l_+) + 1,$$
because of the extra parameter $r$.
The deformation of homotopies is said to be regular if the linearized operator for \eqref{lagSeidel:homotopyeq} is surjective for all $(r,u)$ and no bubbling off of either spheres or discs occur for the moduli spaces with dimension ${\leqslant} 1$.
We note that here the monotonicity guarantees the existence of regular deformation of homotopies.
\qed

\subsection{Module property and Seidel element}\label{lagSeidel:moduleproperty}
\begin{prop}\label{lagSeidel:moduleprop}
The map $FH_*(\tilde g)$ is a module map with respect to the half pair of pants product on $FH_*(M, L)$, i.e. for $[\tilde l_-], [\tilde l_0] \in FH_*(M, L)$, we have
$$FH_*(\tilde g)([\tilde l_-] * [\tilde l_0]) = FH_*(\tilde g)([\tilde l_-]) * [\tilde l_0]$$
\end{prop}
{\it Proof:}
Because of the homotopy invariance, we may reparametrize $g$ so that $g_t = id$ for $t \in [0, \frac{1}{2}]$. Consider the half pair of pants product defined by the punctured strip as in Figure \ref{fig:lagFloer:ppdomain}, \S\ref{lagFloer:prod}. Let $(\mathbf H, \mathbf J)$ be a regular pair which pulls back to the ends $e_\pm$ and $e_0$ respectively as $(H_\pm, J_\pm)$ and $(H_0, J_0)$. Then the pair $(\mathbf H^g, \mathbf J^g)$ defined by
$$\mathbf H^g(s, t, x) := \mathbf H(s, t, g_tx) - K(s, t, g_tx) \text{ and } \mathbf J^g(s, t, x) := dg_t^{-1} \circ \mathbf J_{s, t}(x) \circ d g_t$$
pulls back to the ends $e_\pm$ and $e_0$ respectively as $(H^g_\pm, J^g_\pm)$ and $(H_0, J_0)$. 
Let $\tilde l^g_\pm$ and $\tilde l_\pm$ be critical points of the action functionals $a_{H_\pm}$ and $a_{H^g_\pm}$ respectively and $\tilde l_0$ a critical point of $a_{H_0}$.  We then have the isomorphism of moduli spaces as in Proposition \ref{lagSeidel:fhaction}
$$\mathcal M_{\mathbf H^g, \mathbf J^g}(M, L; \tilde l_-, \tilde l_0, \tilde l_+) \cong \mathcal M_{\mathbf H, \mathbf J}(M, L; \tilde l_-^g, \tilde l_0, \tilde l_+^g): u \mapsto u^g,$$
where $u^g(s, t) := g_t \circ u(s, t)$.
The statement then follows.
\qed

From the properties in proposition \ref{lagSeidel:fhaction} and the homotopy invariance, we may define similarly the Seidel element for the Lagrangian Floer homology. Here we have to assume more:
\begin{definition}\label{lagSeidel:identity}
Suppose that $FH_*(M, L)$ is non-zero and has an identity element with respect to the half-pair-of-pants product defined above, which is denoted $\idelt_L$.
Then 
$$F\Psi_{\widetilde g, L} := FH_*(\widetilde g)(\idelt_L) \in FH_*(M, L)$$
is the Seidel element for the class $[\widetilde g] \in \widetilde\pi_1({\rm Ham}(M, \omega), Ham_L(M, \omega))$.
\end{definition}
\noindent
It follows that in this case, for $[\widetilde l] \in FH_*(M, L)$ we have
$$FH_*(\widetilde g)([\widetilde l]) = F\Psi_{\widetilde g, L} * [\widetilde l].$$
The assumption above is satisfied in many cases, for example the diagonal in $M \times M$.

\subsection{Hamiltonian fibrations over a disc}\label{lagSeidel:hamfib}
The unit disc $D^2$ in ${\mathbb C}$ can be parametrized by the upper half plane $\bar {\mathbb{H}}$ compactified by ${\mathbb R}$ and a point at $\infty$. 
\begin{notation}The following notations are only used in this section. They are NOT compatible with the notations for the same objects used elsewhere in this paper. In the parametrization by ${\mathbb{H}}$, let
$$D^2_{\pm} = \{z \in \bar {\mathbb{H}} | \pm(|z| - 1) {\geqslant} 0\}.$$
\end{notation}
The two half discs can be identified by the map:
$$D^2_+ \to D^2_- : z \mapsto {\bar z}^{-1} \text{ or } re^{i\theta} \mapsto r^{-1}e^{i\theta}.$$

We consider the fibration over $D^2$ defined from an element $g \in {\mathcal{P}}_L{\rm Ham}(M,\omega)$:
$$P_{g} = M \times D^2_{+} \sqcup M \times {D^2_{-}} / \sim : (x, e^{i\pi t}) \sim (g_{t}(x), e^{i\pi t}) \text{ for } t \in [0, 1].$$
Let $\pi : P_g \to D^2$ denote the projection. We note that along the $S^1$-boundary, we have the restricted bundle $N$ that is obtained as the union of  the copies of $L$ in each fiber; it is a Lagrangian submanifold of $P$. Note that $N \simeq L \times S^1$ over $S^1$ if the restriction of $g_1$ to $L$ is diffeotopic to the identity.
A choice of the lifting $\widetilde g \in \widetilde {\mathcal{P}}_L{\rm Ham}(M,\omega)$ amounts again to the choice of a section class $\sigma_{\widetilde g}$ in $\pi_2(P_{g}, N)$ as follows. 
For $x \in L$, consider $\widetilde l = [x, x] \in \widetilde {\mathcal{P}}_L$ and let $\widetilde l^{\widetilde g} = [g_t(x), w]$ with
$$w : D^2_+ \to M : w(e^{i\pi t}) = g_t(x).$$
Via the identification of $D^2_\pm$, we write 
$$w_- : D^2_- \to M: w_-(z) = w({\bar z}^{-1}),$$
in particular, $w_-(e^{i\pi t}) = g_t(x)$ as well. Now the following section in $P_g$ represents $\sigma_{\widetilde g}$:
$$\{x\} \sqcup \{w_-\} / \sim : (x, e^{i\pi t}) \sim (g_{t}(x), e^{i\pi t}) \text{ for } t \in [0,1],$$
where, for example, $\{w_-\}$ denotes the graph of the map $w_-$.
\begin{definition}\label{lagSeidel:vertmaslov}
Let the smooth map $u: D^2 \to P_g$ represent $B \in \pi_2(P_g, N)$. The \emph{vertical Maslov index} of $B$, denoted $\mu^v(B)$ is the Maslov index of the bundle pair $(u^*T^vP_g, (\partial u)^*T^vN)$, where $T^v = \ker d\pi$ denotes the respective vertical tangent bundles.
\end{definition}
\noindent
It's not hard to show that the above is well defined and not dependent on the choice of a smooth map $u$. We then have the following
\begin{prop}\label{lagSeidel:maslovdeg}
$\mu(\widetilde g) = \mu^v([\sigma_{\widetilde g}])$.
\end{prop}
{\it Proof:} 
The trivial Lagrangian path $T^v_xN = T_xL$ over $D^2_+$ is isotopied to $G_t(x){\mathbb R}^n$ over $D^2_-$, via any trivialization chosen for $u^*T^vP_{g}$. The proposition follows from the definitions. 
\qed

\begin{remark}
\rm{Since we will not need it in this article, we leave to the reader to check  that the definition of the action of the paths in ${\mathcal{P}}_L{\rm Ham}(M,\omega)$ on the relative Floer homology can be interpreted in a geometric way on this bundle over the 2-disc, roughly as  the absolute Seidel morphism was interpreted in Lalonde-McDuff-Polterovich \cite{LalondeMcDuffPolterovich} as a map from the quantum homology of the fiber at the north pole to the quantum homology of the fiber at the south pole in a fibration over the 2-sphere. 

For this purpose, one considers the fibers $(M_{\pm 1}, L_{\pm 1}) = \pi^{-1}(\pm 1)$. The natural map from $FH_*(M_1, L_1)$ to $FH_*(M_{-1}, L_{-1})$ can be defined by the pearl complex (i.e. linear clusters). Namely, one flows inside $M_1$ from a critical point of the Morse function on $L_1$ in a linear cluster until that cluster reaches a pseudo-holomorphic section $\sigma$ of $P$ with boundary on $N$; then flows along a linear cluster in the fiber $M_{-1}$, starting from the point $\sigma \cap M_{-1} \in L_{-1}$, until it reaches some critical point of the Morse function on $L_{-1} \subset M_{-1}$.
}
\end{remark}

\subsection{Compatibility among the actions}\label{lagSeidel:compatiblity}

We start by noting the obvious inclusion:
$$\Omega_0{\rm Ham}(M, \omega) \subset {\mathcal{P}}_L{\rm Ham}(M, \omega),$$
where $\Omega_0{\rm Ham}(M, \omega)$ denotes the group of smooth loops in ${\rm Ham}(M, \omega)$ based at the identity. 
Recall that in \cite{Seidel}, the covering $\widetilde \Omega_0{\rm Ham}(M, \omega)$ is defined as following (cf. proposition \ref{lagSeidel:lifting}):
$$\widetilde \Omega_0{\rm Ham}(M, \omega) := \{\left.(g, \widetilde g) \in \Omega_0{\rm Ham}(M, \omega) \times {\rm Homeo}(\widetilde \Omega M) \right| \widetilde g \text{ lifts the action of } g\}.$$
\begin{lemma}\label{lagSeidel:groupinclusion}
We have the inclusion of groups
$$\widetilde\Omega_0{\rm Ham}(M, \omega) \subset \widetilde{\mathcal{P}}_L{\rm Ham}(M, \omega),$$
extending the inclusion $\Gamma_\omega \xrightarrow i \Gamma_L$ in \S\ref{lagFloer:novikov}. 
\end{lemma}
{\it Proof:} We show that
$$\widetilde\Omega_0{\rm Ham}(M, \omega) \subset \widetilde{\mathcal{P}}^0_L{\rm Ham}(M, \omega),$$
where
$$\widetilde {\mathcal{P}}^0_L{\rm Ham}(M, \omega) := \left\{\left.(g, \widetilde g) \in \widetilde {\mathcal{P}}_L{\rm Ham}(M, \omega) \right| g_1|_L = {\rm id} \in {\rm Diff}(L)\right\}.$$

Let 
$$\Omega_L M = \Omega M \cap {\mathcal{P}}_L M$$
be the space of loops in $M$ starting at points in $L$. Then an element of $\Omega_0{\rm Ham}(M, \omega)$ or ${\mathcal{P}}_L{\rm Ham}(M, \omega)$ is determined by how it acts on $\Omega_L M$. This fact gives a definition of the inclusion $\Omega_0{\rm Ham}(M, \omega) \hookrightarrow  {\mathcal{P}}_L{\rm Ham}(M, \omega)$. 

Let $\pi : \widetilde \Omega M \to \Omega M$ and $\pi_L : \widetilde {\mathcal{P}}_L M \to {\mathcal{P}}_L M$ be the covering projections. Consider
$$\widetilde \Omega_L M := \pi^{-1}(\Omega_L M) \text{ and } \widetilde {\mathcal{P}}_L^0 M := \pi_L^{-1}(\Omega_LM).$$
Then by definition, we have
$$\widetilde \Omega_LM = \{(l, w_\Omega) | l \in \Omega_LM \text{ and } w_\Omega : (D^2, S^1) \to (M, l)\} / \sim_\Omega$$
$$\widetilde {\mathcal{P}}_L^0M = \{(l, w_{\mathcal{P}}) | l \in \Omega_LM \text{ and } w_{\mathcal{P}} : (D^2_+; \partial_+, \partial_0) \to (M; l, L)\} / \sim_{\mathcal{P}}$$
where 
$$w_\Omega \sim_\Omega w'_\Omega \iff I_\omega(v_\Omega) = I_c(v_\Omega) = 0 \text{ for } v_\Omega = w_\Omega \#(-w'_\Omega)$$
$$w_{\mathcal{P}} \sim_{\mathcal{P}} w'_{\mathcal{P}} \iff I_\omega(v_{\mathcal{P}}) = I_\mu(v_{\mathcal{P}}) = 0 \text{ for } v_{\mathcal{P}} = w_{\mathcal{P}} \#(-w'_{\mathcal{P}})$$
Let's choose and fix a smooth map 
$$\iota: (D^2_+; \partial_+, \partial_0) \to (D^2; S^1, \{1\})$$
which contracts $\partial_0$ to $\{1\}$ and is an isomorphism otherwise.
We have for $w_\Omega$
$$\widetilde w_\Omega := w_\Omega \circ \iota: (D^2_+; \partial_+, \partial_0) \to (M; l, L),$$
as well as
$$w_\Omega \sim_\Omega w'_\Omega \iff \widetilde w_\Omega \sim_{\mathcal{P}} \widetilde w'_\Omega.$$
The ``$\Rightarrow$'' above is obvious. The ``$\Leftarrow$'' is because that $I_\mu = 2 I_c$ on the maps of the form $\widetilde w_\Omega \#(- \widetilde w'_\Omega)$.
In particular, $\iota$ induces an inclusion $\iota_*: \widetilde \Omega_LM \to \widetilde {\mathcal{P}}_L^0 M$.

On the other hand, for $w_{\mathcal{P}}$ as above, we define $\partial_0 w_{\mathcal{P}}$ by
$$w_{\mathcal{P}}|_{\partial_0}: ([-1, 1], \{\pm 1\}) \to ([-1, 1]/\{\pm 1\}, \{[1]\}) \xrightarrow{\partial_0 w_{\mathcal{P}}} (L, l(0)= l(1)).$$
We then see that $\partial_0 : w_{\mathcal{P}} \mapsto \partial_0 w_{\mathcal{P}}$ defines a map
$$\partial_0^* : \widetilde {\mathcal{P}}_L^0M \to \pi_1(L) / K$$
where $K$ is the image of $\ker I_\omega \cap \ker I_\mu$ under the map $\pi_2(M, L) \to \pi_1(L)$, and there is the exact sequence
$$0 \to \widetilde \Omega_L M \xrightarrow {\iota_*} \widetilde {\mathcal{P}}_L^0 M \xrightarrow {\partial_0^*} \pi_1(L) / K.$$
It follows that $\widetilde {\mathcal{P}}_L^0M$ is a disjoint union of copies of $\widetilde \Omega_L M$.

Now an element in $\widetilde \Omega_0{\rm Ham}(M, \omega)$ is determined by its action on $\widetilde \Omega_LM$ and one in $\widetilde {\mathcal{P}}^0_L{\rm Ham}(M, \omega)$ by its action on $\widetilde {\mathcal{P}}_L^0M$. 
It follows that $\widetilde \Omega_0{\rm Ham}(M, \omega)$ is the subgroup of $\widetilde {\mathcal{P}}^0_L{\rm Ham}(M, \omega)$ preserving each copy of $\widetilde \Omega_LM$ in $\widetilde {\mathcal{P}}_L^0M$. The rest of the statement is obvious.
\qed

\begin{remark}\label{lagSeidel:homotopysequence}
\rm{From the lemma, we obtain the exact sequence described in \eqref{intro:commutdiag}:
$$\widetilde \pi_1{\rm Ham}(M, \omega) \to \widetilde \pi_1({\rm Ham}(M, \omega), {\rm Ham}_L(M, \omega)) \to \widetilde \pi_0{\rm Ham}_L(M, \omega) \to 0,$$
where the third term is the quotient. From the extension sequences of the first two groups and from the triviality of $\pi_0{\rm Ham}(M, \omega)$, we have the following extension sequence:
$$0 \to \Gamma' \to \widetilde \pi_0{\rm Ham}_L(M, \omega) \to \pi_0{\rm Ham}_L(M, \omega) \to 0,$$
where $\Gamma'$ is a quotient of $\Gamma_L$.
}
\end{remark}

\begin{figure}[!h]
\includegraphics[width=0.9\textwidth]{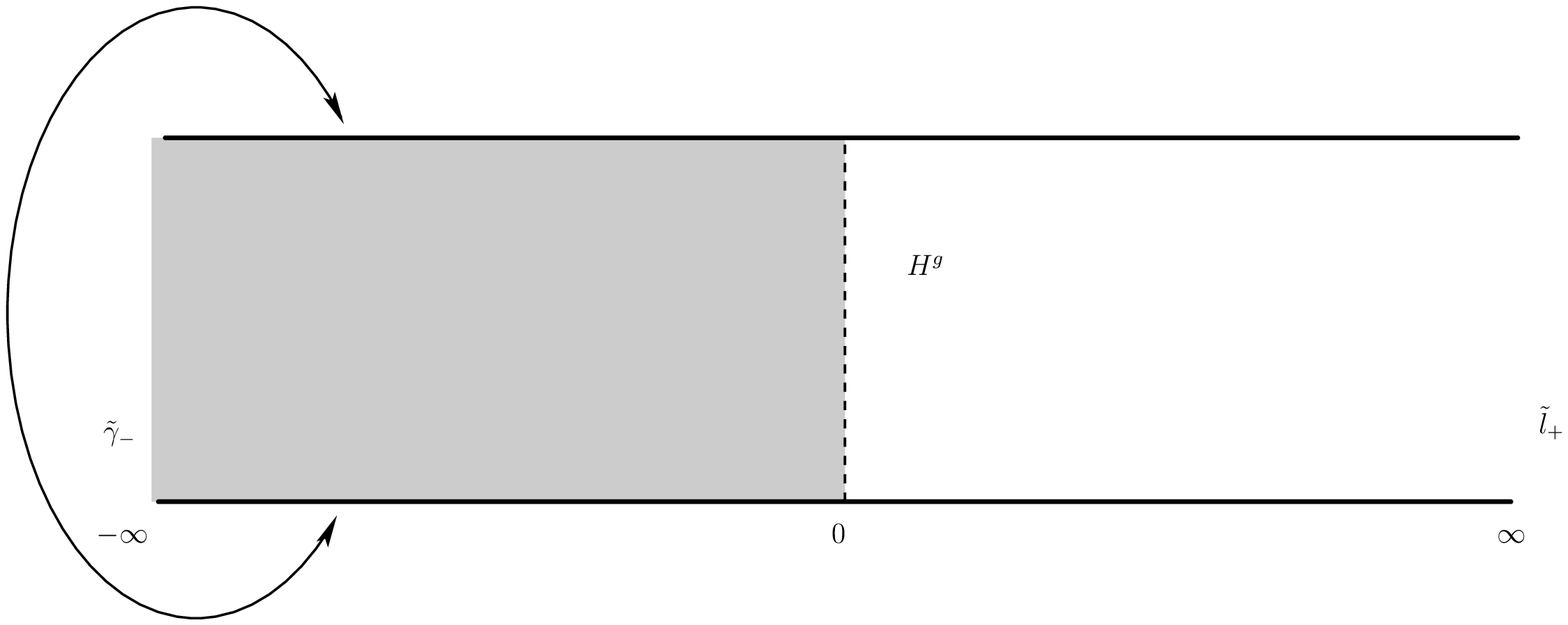}
\caption{}
\label{fig:lagSeidel:schwarzdomain}
\end{figure}

\begin{theorem}\label{lagSeidel:albersseidel}
Let $[\widetilde \gamma] \in FH_*(M, \omega)$ and $\widetilde g \in \widetilde \Omega_0{\rm Ham}(M, \omega) \subset \widetilde {\mathcal{P}}_L{\rm Ham}(M, \omega)$. Then we have \footnote{The same holds for Seidel elements in quantum homology, for which one apply the construction in \cite{McDuff} for the Hamiltonian Seidel elements, while the relative version is obtained from a similar construction in the fibration over a disc. Correspondingly, one needs to consider $H_2^S(M;\mathbb R)$ and $H_2^S(M, L; \mathbb R)$ instead of $H_2^S(M)$ and $H_2^S(M, L)$ when defining the Novikov rings.}
$$\mathscr A(FH_*(\widetilde g)([\widetilde \gamma])) = FH_*(\widetilde g)\mathscr A([\widetilde \gamma]).$$
\end{theorem}
{\it Proof:}
Recall that the description of the ``chimney domain'' in \cite{AbbondandoloSchwarz} as ${\mathbb R} \times [0,1] / \sim$ where $(s, 0) \sim (s,1)$ when $s {\leqslant} 0$, and the conformal structure at $(0,0)$ is given by $\sqrt z$. The domain is depicted in figure \ref{fig:lagSeidel:schwarzdomain}, where the shaded left half of the strip has its two boundaries glued together forming a half infinite cylinder. Let $(H, \mathbf J)$ be a regular pair for both $FH_*(M)$ and $FH_*(M, L)$ so that $(H^g, \mathbf J^g)$ is also regular for both of the theories. 
We then consider an equation similar to \eqref{lagFloer:floweq}:
\begin{equation}\label{lagSeidel:chimneyeq}
\left\{\begin{matrix}
\frac{\partial u}{\partial s} + J^g_t(u)
 \left(\frac{\partial u}{\partial t} - X_{H^g_t}(u)\right) = 0 & \text{ for all }
 (s, t) \in {\mathbb R} \times [0,1], \\
u(s, 0) = u(s, 1) & \text{ for } s \leqslant 0 \\
u|_{[0, \infty) \times \{0, 1\}} \subset L
\end{matrix}\right.
\end{equation}
Then $\mathscr A(\tilde \gamma)$ is defined by counting $0$-dimensional moduli spaces of solutions to \eqref{lagSeidel:chimneyeq}.

We note that $FH_*(\tilde g)([\tilde \gamma])$ is represented by $\tilde \gamma^g$ while $FH_*(\tilde g)([\tilde l])$ by $\tilde l^g$. As in Proposition \ref{lagSeidel:fhaction}, there is a bijection of moduli spaces of solutions to \eqref{lagSeidel:chimneyeq} with the pair $(H^g, \mathbf J^g)$ and that with the pair $(H, \mathbf J)$, given by
$$u \mapsto u^g, \text{ where } u^g (s, t) := g_t \circ u(s, t) : {\mathbb R} \times [0,1] / \sim \to M.$$
This can be shown by directly computing the corresponding equation \eqref{lagSeidel:chimneyeq}. 

It follows that $\mathscr A$ and $FH_*(\tilde g)$ commute on chain level. The independence of choices as well as of $\tilde g$ in the same homotopy class is established similarly as in the case for the Seidel maps.
\qed

\begin{corollary}\label{lagSeidel:3actions}
Assume that the identity exists for $FH_*(M, L)$, then the Seidel element $F\Psi_{\widetilde g, L} \in FH_*(M, L)$ is defined. Let $F\Psi_{\widetilde g}$ denote the Seidel element in $FH_*(M)$ then we have
$$\mathscr A(F\Psi_{\widetilde g}) = F\Psi_{\widetilde g, L}.$$
\end{corollary}
{\it Proof:}
Recall that Proposition \ref{lagFloer:albersidentity} states that $\mathscr A(\idelt) = \idelt_L$ where $\idelt$ and $\idelt_L$ are the respective identity element in $FH_*(M)$ and $FH_*(M, L)$. 
Now replace $[\widetilde \gamma]$ in proposition \ref{lagSeidel:albersseidel} by $\idelt \in FH_*(M)$ and we obtain the proposition.
\qed

\begin{remark}\label{lagSeidel:commutdiag}
\rm{
The above corollary completes the commutative diagram \eqref{intro:commutdiag}, where the maps $\Psi$ and $\Psi_L$ are defined respectively as
$$\Psi(\widetilde g) := F\Psi_{\widetilde g} \text{ and } \Psi_L(\widetilde g) := F\Psi_{\widetilde g, L}.$$
} 
\end{remark}

\section{Reversing the sign of the symplectic structure}\label{reversed:novikov}

We consider here the effects of reversing the symplectic structure on $M$, i.e. the relations between the structures defined on $(M, \omega)$ and $(M, -\omega)$.

Fix a Lagrangian submanifold $L \subset (M, \omega)$.
Let $\omega' = -\omega$ and $c_1'$ and $\mu'$ denote respectively the Chern class and Maslov class for the reversed symplectic structure and $L$, then we obviously have
$$I_{\omega'} = - I_\omega, I_{c'} = -I_c \text{ and } I_{\mu'} = -I_\mu.$$
Correspondingly, we have the Novikov rings $\Lambda_{\omega'}$ and $\Lambda_{L'}$.
Let $\tau : \pi_2(M) \to \pi_2(M)$ and $\pi_2(M, L) \to \pi_2(M, L)$ be the respective involution induced by reversing the signs, then it induces involutions $\tau$ of the groups $\Gamma_\omega$ and $\Gamma_L$ as well as isomorphisms of the Novikov rings as graded rings:
\begin{equation}\label{reversed:novikoviso} 
\tau: \Lambda_{\omega} \to \Lambda_{\omega'} \text{ and } \tau: \Lambda_{L} \to \Lambda_{L'} : 
a_B e^B \mapsto (-1)^{\frac{1}{2}\deg e^{B}} a_{B} e^{\tau (B)}.
\end{equation}
Under our assumption, we see that $\deg e^B$ is always even for either of the two Novikov rings, and thus the above is an isomorphism over ${\mathbb R}$.

\subsection{Quantum ring structure on $QH_*(M)$}\label{reversed:quantring}
Let $2m = \dim_{\mathbb R} M$, then the orientation of $(M, \omega')$ is the $(-1)^m$-multiple of that of $(M, \omega)$. 

\begin{lemma}\label{reversed:classicalring}
Let $\pitchfork$ and $\pitchfork'$ denote the intersections products on $H_*(M, \omega)$ and $H_*(M, \omega')$ respectively. Then we have 
$$\tau(\alpha \pitchfork \beta) = \tau(\alpha) \pitchfork' \tau(\beta)$$
where $\alpha, \beta \in H_*(M) = H_*(M, \omega) = H_*(M, \omega')$ and
$$\tau : H_*(M, \omega) \to H_*(M, \omega') : 
\alpha \mapsto (-1)^m \alpha.
$$
\end{lemma}
{\it Proof:} Let $\{\gamma_j\}$ be a base of $H_*(M)$ and $\{\gamma_j^*\}$ it's dual base with respect to the product $\pitchfork$ and $\{\gamma_j^{*'}\}$ that with respect to $\pitchfork'$. Thus we have
$$\gamma_j^{*'} = (-1)^m \gamma_j^*.$$
Let $a, b, c_j$ and $c_j^{*}$ be generic cycles representing $\alpha, \beta, \gamma_j$ and $\gamma_j^{*}$.
The intersection product $\pitchfork$ (respectively $\pitchfork'$) is alternatively written as
$$\alpha \pitchfork \beta = \sum_j {\langle}\alpha, \beta, \gamma_j^*{\rangle}\gamma_j \text{ (respectively } \alpha\pitchfork' \beta = \sum_j {\langle}\alpha, \beta, \gamma_j^{*'}{\rangle}' \gamma_j \text{)},$$
where ${\langle}\alpha, \beta, \gamma_j^*{\rangle}$ is the intersection number of $a\times b\times c_j^*$ with $\triangle$, the minimal diagonal, in $M^3$, oriented by $\omega$. We only need to compare the coefficients in front of the $\gamma_j$'s.

The orientations of the cycle $\triangle$ in $(M, \omega)^3$ and $(M, -\omega)^3$, as well as the orientations of $M^3$ in either case, differ by $(-1)^m$, while the orientations of $\gamma_j^*$ and $\gamma_j^{*'}$ also differ by $(-1)^m$. 
The lemma then follows.
\qed

We consider the effect on $QH_*(M)$. Since
$$QH_*(M, \omega) = H_*(M) \otimes \Lambda_\omega \text{ and } \tau: \Lambda_\omega \to \Lambda_{\omega'},$$
we define naturally the induced map by:
\begin{equation}\label{reversed:quantumcomap}\tau_*: QH_*(M, \omega) \to QH_*(M, \omega') : \alpha\otimes f \mapsto \tau(\alpha) \otimes \tau(f).\end{equation}

\begin{prop}\label{reversed:quantumcoho}
The map $\tau_*$ defined in \eqref{reversed:quantumcomap} is a ring isomorphism of quantum homologies over the isomorphism $\tau$ of the Novikov rings in \eqref{reversed:novikoviso}.
\end{prop}
{\it Proof:}
The quantum intersection product on $QH_*(M, \omega)$ (resp. $QH_*(M, \omega')$) is denoted $*$ (resp. $*'$). Choose and fix a base $\{\gamma_j\}$ of $H_*(M)$ and denote $\{\gamma_j^*\}$ (resp. $\{\gamma_j^{*'}\}$) the dual base with respect to the intersection product on $(M, \omega)$ (resp. $(M, -\omega$)). Then for $\alpha, \beta \in H_*(M)$:
$$\alpha * \beta = \sum_{j, B} {\langle}\alpha, \beta, \gamma^*_j{\rangle}_B \gamma_j e^{-B} \text{ and } \alpha *' \beta = \sum_{j, B} {\langle}\alpha, \beta, \gamma_j^{*'}{\rangle}'_B \gamma_j e^{-B}.$$
We need to check that $\tau_*\{(\alpha e^A) * (\beta e^B)\} = \tau_*(\alpha e^A) *' \tau_*(\beta e^B)$, where we dropped the $\otimes$ in the expressions.
It follows from lemma \ref{reversed:qhcoefcompare} below, which compares the coefficients of the two quantum intersection products. Given the lemma, we have
\begin{equation*}
\begin{split}
LHS = & \tau_*\{(\alpha * \beta) e^{A+B}\} = \tau_*\left\{\sum_{j, C} {\langle}\alpha, \beta, \gamma_j^*{\rangle}_C \gamma_j e^{A+B-C}\right\} \\
= & \sum_{j, C} (-1)^{m + I_c(A+B-C)} {\langle}\alpha, \beta, \gamma_j^*{\rangle}_C \gamma_j e^{\tau(A+B-C)}\\
\left._\lozenge\right.= & \sum_{j, C} (-1)^{I_c(A+B)}{\langle}\alpha, \beta, \gamma_j^{*'}{\rangle}'_{\tau(C)} \gamma_j e^{\tau(A+B-C)}  \\
= & \sum_{j, C} {\langle}(-1)^m \alpha, (-1)^m \beta, \gamma_j^{*'}{\rangle}'_{\tau(C)} \gamma_j e^{-\tau(C)} \left\{(-1)^{I_c(A)}e^{\tau(A)}\right\} \left\{(-1)^{I_c(B)}e^{\tau(B)}\right\} \\
= RHS &
\end{split}
\end{equation*}
where $\lozenge$ is lemma \ref{reversed:qhcoefcompare}. 
\qed
\begin{lemma}\label{reversed:qhcoefcompare}
For all $B \in \Gamma_\omega$ and $j$, we have ${\langle}\alpha, \beta, \gamma_j^*{\rangle}_B = (-1)^{m + I_c(B)}{\langle}\alpha, \beta, \gamma_j^{*'}{\rangle}'_{\tau(B)}$.
\end{lemma}
{\it Proof:}
We recall first the definition of the triple intersection ${\langle}\alpha, \beta, \gamma_j^*{\rangle}_B$. Consider the moduli space $\bar{\mathcal{M}}_{0,3}(M, \omega, J; B)$ of $J$-holomophic spheres in $M$ with $3$-marked points, representing $B \in \Gamma_\omega$. The marked points fixes the parametrization of the principle components in the domain and we assume that they corresponds to $0, 1$ and $\infty$ (in that order) respectively. Let $ev$ denote the evaluation map
$$ev: \bar{\mathcal{M}}_{0,3}(M, \omega, J; B) \to M^3.$$
Choose and fix generic cycles $a, b$ and $c^*$ in $M$ representing the classes $\alpha, \beta$ and $\gamma_j^*$, so that $ev$ is transversal to $a \times b \times c^*$. Then the triple intersection is defined to be the cardinality of the following intersection when the resulting dimension is $0$:
$${\langle}\alpha, \beta, \gamma_j^*{\rangle}_B = ev_*([\bar{\mathcal{M}}]) \pitchfork a\times b \times c^*.$$

Let $\rho: {\mathbb{CP}}^1 \to {\mathbb{CP}}^1$ denote the standard complex conjugation on ${\mathbb{CP}}^1 = S^2$, in particular, it fixes the $3$ marked points $0, 1$ and $\infty$. We note that $u : S^2 \to (M, \omega, J)$ is $J$-holomorphic and represents $B \in \Gamma_\omega$ iff $\rho(u) : S^2 \xrightarrow{\rho} S^2 \to (M, \omega', J')$ is $J'$-holomorphic and represents $\tau(B) \in \Gamma_{\omega'}$. We can in fact establish an explicit identification of the moduli spaces:
$$\rho :  \bar{\mathcal{M}}_{0,3}(M, \omega, J; B) \to  \bar{\mathcal{M}}_{0,3}(M, \omega', J'; \tau(B)) : u \mapsto \rho(u),$$
where slightly more care is taken in case of the nodal domains. Furthermore, the evaluation maps coincide, i.e.
$$ev = ev' \circ \bar \rho, \text{ where } ev' :  \bar{\mathcal{M}}_{0,3}(M, \omega', J'; \tau(B)) \to M^3.$$
It then follows that the two triple intersections coincide upto a sign.

The sign comes from two sources, the manifold $M$ and the moduli spaces ${\mathcal{M}}$. The orientation of $(M, -\omega)$ implies that
$$\gamma_j^{*'} = (-1)^m \gamma_j^* \text{ and } \pitchfork' = (-1)^{3m} \pitchfork.$$
It follows that the overall sign only comes from ${\mathcal{M}}$ and the identification $\rho$. 
We check this sign in the following. Fix $u \in {\mathcal{M}}_{0,3}(M, \omega, J; B)$, then the tangent space at $u$ is given by the linearized operator $D\bar\partial_J$:
$$D\bar\partial_J(\xi)(z) = \nabla\xi(z) + J(u(z)) \circ \nabla \xi(z) \circ j_z + l.o.t., \text{ for } \xi \in \Omega^0(u^*TM), z \in S^2,$$
where $\nabla$ is the induced connection on $u^*TM$ from a Hermitian connection on $TM$, compatible with $(\omega, J)$. The operator $D\bar\partial_J$ can be homotopied through Fredholm operators to the standard $\bar\partial$ operator on the holomorphic vector bundle $u^*T^{1,0}_JM$. Under this homotopy, we obtain an identification of the solution space to $H^0({\mathbb{CP}}^1, u^*T^{1,0}_{J} M)$.
The orientation of ${\mathcal{M}}_{0,3}(M, \omega, J; B)$ at $u$ is then defined by the canonical orientation of the complex vector space $H^0({\mathbb{CP}}^1, u^*T^{1,0}_{J} M)$.

For $v := \rho(u) \in {\mathcal{M}}_{0,3}(M, \omega', J'; \tau(B))$, we have similarly the linearized operator
$$D\bar\partial_{J'}(\zeta)(z) = \nabla' \zeta(z) + J'(v(z)) \circ \nabla' \zeta(z) \circ j_z + l.o.t, \text{ for } \zeta \in \Omega^0(v^*TM),$$
where $\nabla'$ is the induced connection on $v^*TM$ from the same Hermitian connection on $TM$. The following in fact holds:
$$D\bar\partial_{J'} = \rho^*D\bar\partial_J,$$
and thus the homotopy to $\bar\partial$ is pulled back via $\rho$. The orientation of the moduli space ${\mathcal{M}}_{0,3}(M, \omega', J'; \tau(B))$ at $v$ is thus defined by the canonical orientation of the complex vector space $H^0({\mathbb{CP}}^1, v^*T^{1,0}_{J'}(M))$.

The tangent map $d\rho$ at $u$ is:
$$d\rho : \xi \mapsto \rho^*\xi \text{ where } (\rho^*\xi)(z) = \xi(\rho(z)),$$
which induces the following indentification as real vector spaces:
$$d\rho : H^0({\mathbb{CP}}^1, E) \to H^0({\mathbb{CP}}^1, \rho^*\bar E),$$
where $E = u^*T^{0,1}_J M$ is a rank $n$ holomorphic vector bundle over ${\mathbb{CP}}^1$. We check that $d\rho$ is complex anti-linear by evaluating $\xi$ and $d\rho(\xi)$ at respective points in ${\mathbb{CP}}^1$. Since the fibers $E_z$ and $(\rho^*\bar E)_{\rho(z)}$ are identical with opposite complex structures, we have for $\lambda \in {\mathbb C}$:
$$(\rho^*(\lambda \xi)_E)(z) = (\lambda\xi)_E(\rho(z)) = (\bar\lambda\xi)_{\bar E}(\rho(z)) = (\bar\lambda \rho^*\xi)_{\bar E}(z).$$
It follows that the orientation of the map $\rho$ is given by $(-1)^{\dim {\mathcal{M}}_{0,3}(M, B)} = (-1)^{m + I_c(B)}$ and the lemma follows.
\qed

\subsection{Seidel elements in $QH_*(M)$}\label{reversed:seidelelts}
The group ${\rm Ham}(M, \omega)$ is naturally a subgroup of ${\rm Diff}(M)$. Suppose that $\{H_t\}_{t \in [0,1]}$ generates $g_{t \in [0,1]} \in {\rm Ham}(M, \omega)$, then regarded as an element in ${\rm Ham}(M, -\omega)$, it is alternatively generated by $\{-H_t\}_{t \in [0,1]}$.
We see that
$${\rm Ham}(M, \omega) = {\rm Ham}(M, -\omega) \subset {\rm Diff}(M).$$

We define a \emph{reversion map} $\tau$ on the group of loops:
\begin{equation}\label{reversed:loopgroup}\tau : \Omega_0{\rm Ham}(M, \omega) \to \Omega_0{\rm Ham}(M, -\omega) : g := \{g_t\} \mapsto g^- := \{g_{1-t}\}.
\end{equation}
The following lemma is obvious:
\begin{lemma}\label{reversed:genham}
Suppose $K = K_t$ generates the loop $g \in \Omega_0{\rm Ham}(M, \omega)$, then $\underline K := \{K_{1-t}\}$ generates the loop $g^- \in \Omega_0(M, -\omega)$.
\qed
\end{lemma}
\noindent
We note that the loop $g^-$ is homotopic to $g^{-1}$ in ${\rm Ham}(M, \pm \omega)$, viewed as identical subgroup in ${\rm Diff}(M)$. 

The reversion map $\tau$ on $\Omega M$ can be extended to $\widetilde \Omega M$ via:
$$\tau([\gamma, v]) = [\tau(\gamma), \tau(v)] \text{ where } \tau(v) : D^2 \to M : z \mapsto v(\bar z).$$
Then $\tau : \widetilde \Omega_0{\rm Ham}(M, \omega) \to \widetilde \Omega_0{\rm Ham}(M, -\omega)$ is defined by
$$\tau(g, \widetilde g) (\tau(\widetilde \gamma)) = \tau \circ (g, \widetilde g) \circ \tau (\widetilde \gamma) \text{ for } \widetilde \gamma \in \widetilde \Omega M.$$

In the following, we will use the description of the Seidel element $\Psi_{[g]}$ for $[g] \in \pi_1{\rm Ham}(M, \omega)$ in terms of Gromov-Witten invariants in the Hamiltonian fibration $P_{[g]} \to S^2$ defined from $[g]$ as in Lalonde-McDuff-Polterovich \cite{LalondeMcDuffPolterovich}.
\begin{prop}\label{reversed:hamSeidel}
$$\Psi_{\tau([g])} = \tau_*(\Psi_{[g]}) \in QH_*(M, -\omega).$$
\end{prop}
\noindent
{\it Proof:}
Let $(g, \widetilde g) \in \widetilde{\Omega}_0{\rm Ham}(M, \omega)$ and $\Psi_{\widetilde g}$ the corresponding Seidel element.
Let $P_g \xrightarrow \pi S^2$ be the fibration defined by $g$:
$$P_g = D^2_1\times M \cup_g D^2_2 \times M, \text{ where } D^2_1 \times M \ni (e^{2\pi it}, x) \sim (e^{2\pi it}, g_t(x)) \in D^2_2 \times M,$$
and $\kappa$ the coupling form on $P_g$, extending $\omega$ on the fibers, then, for appropriate $\varepsilon > 0$, $\omega_g = \pi^*\omega_0 + \varepsilon \kappa$ is a symplectic form on $P_g$ where $\omega_0$ is the standard symplectic form on $S^2$ inducing the positive orientation on $D^2_1$ (thus negative orientation on $D^2_2$). Then $\Psi_{\widetilde g}$ is defined by looking at the section classes in $P_g$.

The corresponding bundle $P_{g^-}$ can be defined similarly. We give an alternative construction below. Let $r : D^2 \to D^2$ be the standard conjugation as the unit disc in ${\mathbb C}$, and use the same letter $r$ to denote the induced conjugation map on $S^2 = D^2_1 \cup_\partial D^2_2$. Then $$P_{g^-} = r^*P_g, \kappa^- = -r^*\kappa, P_{g} = r^*P_{g^-}, \kappa = - r^*\kappa^- \text{ and } r^*\omega_0 = -\omega_0,$$ 
where, of course, $r$ is also used to denote the pull-back maps between the Hamiltonian fibrations in the above. The symplectic form on $P_{g^-}$ is then
$$\omega_{g^-} = \pi^*\omega_0 + \varepsilon\kappa^- \Rightarrow r^*\omega_{g^-} = - \omega_{g},$$
i.e. $r : (P_{g^-}, \omega_{g^-}) \xrightarrow \simeq (P_{g}, -\omega_{g})$ symplectically.
The two sides will be used interchangeably.

Let $\sigma_0 \in H_2(P_g)$ be the standard reference section class (cf. \cite{McDuff} Lemma 3.2), for example, when $c_1(TM)$ and $[\omega]$ are not proportional on spherical classes in $M$, we ask
$$c_1^v(\sigma_0) = \kappa(\sigma_0) = 0.$$
Then we have $\Psi_{\widetilde g} = e^{\sigma_{\widetilde g}} \Psi_{[g]}$ for some $\sigma_{\widetilde g} \in H_2(M; {\mathbb R})$ and
$$\Psi_{[g]} = \sum_{B, j} {\langle}[M], [M], \iota_*(\gamma_j^*){\rangle}_{\sigma_0 + \iota_*(B)} \gamma_j e^{-B} \in QH_*(M, \omega),$$
where $\iota : M \to P_g$ is the inclusion of a fiber, $B \in H_2(M; {\mathbb R})$ so that $\sigma_0 + \iota_*(B) \in H_2(P_g; {\mathbb{Z}})$ is represented by a section and $\{\gamma_j\}$, $\{\gamma_j^*\}$ are dual bases of $H_*(M)$ under $\pitchfork$ and ${\langle}\ldots{\rangle}_{\ldots}$ denotes the Gromov-Witten invariants in $(P_g, \omega_g)$. It follows that
$$\tau_*(\Psi_{[g]}) = \sum_{B, j} (-1)^{m + I_c(B)}{\langle}[M], [M], \iota_*(\gamma_j^*){\rangle}_{\sigma_0 + \iota_*(B)} \gamma_j e^{-\tau(B)} \in QH_*(M, -\omega).$$

Let $\sigma$ be a section class in $(P_g, \omega_g)$, i.e. $\pi_*\sigma = [(S^2, \omega_0)]$ in the natural orientation, then $\tau(\sigma) := -\sigma$ is a section class in $(P_g, -\omega_g)$, because $\pi_*(-\sigma) = [(S^2, -\omega_0)]$. On the other hand, $\sigma_0^- := \tau(\sigma_0)$ is a standard reference section class as well.
We may write down the Seidel element for $\tau(\widetilde g)$ in $QH_*(M, -\omega)$ as $\Psi_{\tau(\widetilde g)} = e^{\tau(\sigma_{\widetilde g})} \Psi_{\tau([g])}$, where:
\begin{equation*}
\Psi_{\tau([g])} = \sum_{B, j} {\langle}[M], [M], \iota_*(\gamma_j^{*'}){\rangle}'_{\sigma_0^- + \iota_*(\tau(B))} \gamma_j e^{-\tau(B)},
\end{equation*}
and we have to show that:
$${\langle}[M], [M], \iota_*(\gamma_j^{*'}){\rangle}'_{\sigma_0^- + \iota_*(\tau(B))} =  (-1)^{m +I_c(B)}{\langle}[M], [M], \iota_*(\gamma_j^*){\rangle}_{\sigma_0 + \iota_*(B)}.$$

The dimension of the relevant moduli spaces is $m + 1+c_1(TP)(\sigma_0 + \iota_*(B))$. Let $[P_g, \omega_g]$ be the fundamental class of $P_g$ with orientation given by $\omega_g$, then 
$$[P_g, \omega_g] = (-1)^{m+1} [P_g, -\omega_g] \text{ and }$$
$$[M, \omega] = (-1)^m [M, -\omega].$$
We have also $\gamma_j^{*'} = (-1)^m \gamma_j^*$. It follows from the same argument as in lemma \ref{reversed:qhcoefcompare} that the overall sign for the Gromov-Witten invariants is given by
$$(-1)^{3m+3 + 3m + m+1 + c_1(TP)(\sigma_0 + \iota_*(B))} = (-1)^{m+ c_1(TS^2)([S^2]) + c_1^v(B)} = (-1)^{m + I_c(B)}.$$
\qed

\section{Reversing operations in Lagrangian Floer homology}\label{reversed:lagFloer}

We define first a \emph{reversion map} on $\widetilde {\mathcal{P}}_LM$, that we denote by $\underline\tau$. Let $\widetilde l = [l, w]$ denote a typical element of $\widetilde {\mathcal{P}}_LM$, i.e.
$$l : ([0,1], \{0,1\}) \to (M, L) \text{ and } w: (D^2_+; \partial_+, \partial_0) \to (M; l, L).$$
Then we define $\underline {\widetilde l} := \underline \tau (\widetilde l)$ by
$$\underline l : [0,1] \to M : t \mapsto l(1-t) \text{ and } \underline w: D^2_+ \to M : z \mapsto w(-\bar z).$$
It's obvious that $\underline\tau$ is an involution, i.e. $\underline\tau^2 = id$. We note that the action of $\pi_2(M, L)$ on $\widetilde {\mathcal{P}}_LM$ as deck transformations are intertwined by $\underline\tau$:
$$\underline {(B\circ \widetilde l)} = \underline \tau([l, w \# B]) = [\underline l, \underline {w \# B}] = [\underline l, \underline w \#\tau(B)] = \tau(B) \circ \underline{\widetilde l}.$$
It follows that $\underline\tau$ defines an involution on $\widetilde {\mathcal{P}}_LM$.

Let now $(H, \mathbf J)$ be a regular pair for defining the Floer homology of $(M, L; \omega)$. We consider the \emph{reversed} pair $(\underline H, \underline{\mathbf J})$:
$$\underline H_t = H_{1-t} \text{ and } \underline J_t = -J_{1-t}.$$
Then it's easy to check that the corresponding action functionals satisfy
$$a_{\underline H} (\underline{\widetilde l}) = a_{H}(\widetilde l).$$
In fact, the involution $\tau$ identifies the metric $(,)_{\mathbf J}$ with $(,)_{\underline{\mathbf J}}$ as well. We show that the Floer homologies are identified by $\tau$. 

The next three lemmas are obvious. 

\begin{lemma}\label{reversed:nondegen}
$l$ is Hamiltonian path of $H$ in $(M, \omega)$ $\iff$ $\underline{l}$ is a Hamiltonian path for $\underline H$ in $(M, \omega' = -\omega)$. Furthermore, $\widetilde l$ is non-degenerate $\iff$ $\underline{\widetilde l}$ is non-degenerate.
\end{lemma}
{\it Proof:} Let $X_t = \omega(H_t)$ be the Hamiltonian vector field of $H_t$, then $\underline X_t = -\omega(H_{1-t}) = -X_{1-t}$ is the Hamiltonian vector field of $\underline H_t$. Let $\phi_t$ and $\underline\phi_t$ be the Hamiltonian isotopies generated by $X_t$ and $\underline X_t$ respectively, then $\underline\phi_t = \phi_{1-t}\circ \phi_1^{-1}$. Thus $\underline\phi_1 = \phi_1^{-1}$ and the lemma follows.
\qed

Next we compute the Conley-Zehnder index.
\begin{lemma}\label{reversed:czindex}
$\mu_H(\widetilde l) = \mu_{\underline H}(\underline{\widetilde l})$.
\end{lemma}
{\it Proof:}
We recall the notations in \S\ref{lagFloer:maslov}. Let $\Phi_z: (T_{w(z)} M, \omega) \to ({\mathbb C}^n, \omega_0)$ be the trivialization of $w^*TM$ so that $\Phi_r(T_{w(r)}L) = {\mathbb R}^n$ for all $r \in [-1, 1]$ and 
$$E_t = \Phi_{e^{i\pi t}} \circ d\phi_t \circ \Phi^{-1}_1 \in Sp({\mathbb C}^n, \omega_0)$$
be the path of symplectic matrices. Then
$$\mu_H(\widetilde l) = \mu(E_t{\mathbb R}^n \oplus {\mathbb R}^n, \triangle)$$
where the symplectic structure on ${\mathbb C}^n \oplus {\mathbb C}^n$ is given by $(\omega_0 \oplus -\omega_0)$. For $\underline{\widetilde l}$, a symplectic trivialization of $\underline w^*TM$ is given by
$$\underline \Phi_z := \Phi_{-\bar z} : (T_{\underline w(z)} M, -\omega) \to ({\mathbb C}^n, -\omega_0),$$
and we have
$$\underline E_t := \underline\Phi_{e^{i\pi t}} \circ d\underline\phi_t \circ \underline\Phi^{-1}_1 = E_{1-t} \circ E_1^{-1}.$$

Now the index we need is
$$\mu_{\underline H}(\underline{\widetilde l}) = \mu(\underline E_t {\mathbb R}^n \oplus {\mathbb R}^n, \triangle) \text{ in } ({\mathbb C}^n \oplus {\mathbb C}^n, -\omega_0 \oplus \omega_0).$$
Comparing with the expression for $\mu_H(\widetilde l)$, this reverses both the symplectic structure and the path of symplectic matrices. The property of the Maslov index of pairs as defined in \cite{RobbinSalamon2} implies the lemma.
\qed

\begin{lemma}\label{reversed:nondegenerate}
The pair $(H, \mathbf J)$ is regular iff $(\underline H, \underline{\mathbf J})$ is regular.
\end{lemma}
{\it Proof:}
It is staightforward to check that the defining equations for the various objects involved in either case are identified by transformations induced from $\tau$. 
\qed

By \S\ref{lagFloer:orient}, the orientations of the trajectories are given by those of the moduli spaces of discs (canonically given by the choice of relative spin structure) and the moduli spaces of capped strips. Here, we discuss first the effect of reversion on the moduli spaces of discs. 
We consider the parametrized disc $D^2$. 
Let $\rho: D^2 \to D^2$ be the complex conjugation on $D^2 \subset {\mathbb C}$. Obviously $u : (D^2, S^1) \to (M, L; \omega, J)$ is a holomorphic disc with boundary on $L$ representing $B \in \pi_2(M, L)$ iff $\rho(u) : (D^2, S^1) \xrightarrow \rho (D^2, S^1) \xrightarrow u (M, L; -\omega, -J)$ is holomorphic and represents $\tau(B) \in \Gamma_{L'}$. Let $\widetilde {\mathcal{M}}(M, L; \omega, J; B)$ denote the moduli space of parametrized $J$-holomorphic discs representing the class $B$.
\begin{lemma}
With the same choice of the relative spin structure of $L$ in $M$, the orientation of the map 
$$\rho :\widetilde {\mathcal{M}}(M, L; \omega, J; B)\to \widetilde {\mathcal{M}}(M, L; -\omega, -J; \tau(B)) : u \mapsto \rho(u)$$
is given by $(-1)^{\frac{1}{2}\deg B}$.
\end{lemma}
{\it Proof:}
Recall that the orientation of the moduli space of discs is given by the identification (cf. \cite{FOOO, BiranCornea}):
$$\ker D\bar\partial_J \simeq \ker({\rm Hol} _J(D^2, S^1; {\mathbb C}^n, {\mathbb R}^n) \times {\rm Hol} _J(S^2; E) \xrightarrow{ev} {\mathbb C}^n).$$
The three items on the right are oriented respectively by the following. The first item is oriented by the choice of the relative spin structure, while independent of the structure $J$. With the choice of the relative spin structure, the second item is oriented by the structure $J$. The last item is oriented by $J$ while independent of the relative spin structure. Under the map $\rho$, the first item is canonically identified, while the rest follows similarly as in lemma \ref{reversed:qhcoefcompare}. Thus, if we fix the relative spin structure of $L$ and reverse $J$, the orientation of the moduli space is changed by $(-1)^{\frac{1}{2}\mu_L(B)}$, which by definition is $(-1)^{-\frac{1}{2}\deg B} = (-1)^{\frac{1}{2}\deg B}$.
\qed

The reversing map $\underline\tau$ on $\widetilde {\mathcal{P}}_LM$ induces the correspondence between the respective caps (cf. \eqref{lagFloer:caps}) via the complex conjugation $\rho$ of $Z_\pm \subset {\mathbb C}$. We see that $u^\pm$ is a solution of the equation \eqref{lagFloer:capeq} for $(\omega, \mathbf J, J, H)$ iff $\underline u^\pm = u^\pm \circ \rho$ is a solution for $(-\omega, \underline{\mathbf J}, \underline J, \underline H)$. It follows that the corresponding moduli spaces of caps are isomorphic via the map
$$ \rho: u \mapsto \rho(u) = \underline u.$$
We assign the orientations for the reversed moduli spaces so that the map $\rho$ preserves the orientations for the preferred basis. The orientations of the reversed caps given by the reversed preferred basis are related by
$$(-1)^{\rho(\widetilde l \# B)} = (-1)^{\rho(\widetilde l) + \frac{1}{2}\deg B}.$$

\begin{prop}\label{reversed:isofloer}$\underline\tau$ induces an isomorphism of Floer homologies, intertwining as well the isomorphism of Novikov rings \eqref{reversed:novikoviso},
$$\tau_* : FH_*(M, L, \omega; H, \mathbf J) \to FH_*(M, L, -\omega; \underline H, \underline{\mathbf J}).$$
\end{prop}
{\it Proof:}
The map $\tau$ induces natural transformation taking the equation \eqref{lagFloer:floweq} for the left side to the one for the right side:
$$v(s, t) := \tau(u)(s, t) = u(s, 1-t) \text{ so that}$$
\begin{equation*}
\begin{split}
& \frac{\partial v}{\partial s} + \underline J_t(v) \left(\frac{\partial v}{\partial t} - X_{\underline H_t}(v)\right) = \\
 = & \left.\frac{\partial u}{\partial s}\right|_{1-t} + J_{1-t}(u(s, 1-t)) \left(\left.\frac{\partial u}{\partial t}\right|_{1-t} - X_{H_{1-t}}(u(s, 1-t))\right) = 0
\end{split}
\end{equation*}

Together with the last three lemmas, we see that the moduli spaces, as well as the compactifications, correspond via $\tau$ and $\underline \tau$. 
We then have an isomorphism at the chain level (where the orientations of the moduli spaces identified by $\tau$ are defined to be the same) and thus the proposition follows. The intertwining of the isomorphism \eqref{reversed:novikoviso} is automatic.
\qed

\end{document}